\pdfoutput=1
\documentclass[11pt]{article}
\usepackage{graphicx, amsfonts, eucal}

\newtheorem{Proposition}{Proposition}
\newtheorem{Theorem}{Theorem}
\newtheorem{Lemma}{Lemma}
\newtheorem{Cor}{Corollary}
\newtheorem{Definition}{Definition}

\newcommand{\proof}{\par \noindent {\bf Proof: }}

\newcommand{\rmk}{\bigskip \par \noindent {\em Remark: }}
\newcommand{\eqr}[1]{(\ref{#1})}

\newcommand{\cA}{\mathcal{A}}                 

\newcommand{\cC}{\mathcal{C}}                 
\newcommand{\cE}{\mathcal{E}}

\newcommand{\cP}{\mathcal{P}}                 
                 
\newcommand{\cS}{\mathcal{S}}

\newcommand{\cV}{\mathcal{V}}

\newcommand{\proj}[1]{\bbP_{\!#1}}
\newcommand{\projh}[1]{\proj{{\rm hor}}}

\newcommand{\R}{{\mathbb R}}

\newcommand{\lp}{\left (}
\newcommand{\rp}{\right )}
\newcommand{\la}{\left \langle}
\newcommand{\ra}{\right \rangle}
\newcommand{\lsb}{\left [}
\newcommand{\rsb}{\right ]}
\newcommand{\lcb}{\left \{}
\newcommand{\rcb}{\right \}}
\newcommand{\dla}{\left \langle \! \left \langle}
\newcommand{\dra}{\right \rangle \! \right \rangle}
\newcommand{\sands}{\qquad \mbox{and} \qquad}
\newcommand{\half}{{\textstyle {1 \over 2}}}

\newcommand{\smallfrac}[2]{{\textstyle {#1 \over #2}}}

\newcommand{\beq}[1]{\begin{equation}\label{#1}}
\newcommand{\eeq}{\end{equation}}
\newcommand{\beqa}{\begin{eqnarray*}}
\newcommand{\eeqa}{\end{eqnarray*}}
\newcommand{\beqan}{\begin{eqnarray}}
\newcommand{\eeqan}{\end{eqnarray}}
\newcommand{\Prop}[1]{\acapon\begin{Proposition}\label{#1}}
\newcommand{\eProp}{\end{Proposition}}
\newcommand{\Lem}[1]{\acapon\begin{Lemma}\label{#1}}
\newcommand{\eLem}{\end{Lemma}}

\newcommand{\prf}{\noindent {\bf Proof: }}
\newcommand{\prfend}{\ \vrule height6pt width6pt depth0pt \medskip}

\newcommand{\rmkend}{ \medskip}

\newcommand{\norm}[1]{\left \Vert #1 \right \Vert}

\newcommand{\bbP}{{\mathbb P}}

\newcommand{\dep}{\smallfrac {d \ }{d \epsilon}}
 
\font\elevenrm=cmr10 at 11pt
\newcommand{\Id}{\hbox{{\rm 1}\kern-3.8pt \elevenrm1}}

\begin{document}
\pagestyle{myheadings}
\thispagestyle{plain}

\newcommand{\bP}{{\Bbb P}}
\newcommand{\bbF}{{\Bbb F}}
\newcommand{\cb}{{\cal B}}
\newcommand{\bn}{{\mathcal B}^n}
\newcommand{\allpaths}{{\cal P}}
\newcommand{\lamu}{{\lambda}}
\newcommand{\lamuf}{{\lambda_{\!f}}}
\newcommand{\lamv}{{\kappa}}
\newcommand{\sfin}{s_{\! f}} 
\newcommand{\tfin}{t_{\! f}}
\newcommand{\xfin}{x_{\! f}}
\newcommand{\yfin}{y_{\! f}}
\newcommand{\zfin}{z_{\! f}}
\newcommand{\umax}{u_{\rm max}}
\newcommand{\tilC}{\widetilde C}
\newcommand{\hatC}{\widehat C}
\newcommand{\tilchi}{\widetilde \chi}
\newcommand{\tilfps}{\widetilde F_{\ps}}
\newcommand{\tilH}{\widetilde H}
\newcommand{\chiH}{\widehat H}
\newcommand{\hatsig}{\widehat \sigma}
\newcommand{\tilQ}{q}
\newcommand{\ps}{\psi_z}
\newcommand{\sgn}{\mbox{sgn}\,}
\newcommand{\pd}[2]{\smallfrac {\partial #1}{\partial #2}}
\newcommand{\depe}[1]{\left . \dep #1 \right |_{\epsilon = 0}}
\newcommand{\vertps}[1]{\mbox{vert}_{\ps}\!#1}
\newcommand{\dddt}{\smallfrac {d^2 \ }{d t^2}}
\newcommand{\tqm}{t_{q, \mu}}
\newcommand{\nqm}{n_{\alpha, p}}
\newcommand{\wqm}{w_{\alpha, p}}
\newcommand{\yqm}{y_{q, \mu}}
\newcommand{\shp}[1]{#1_{{\tiny 1 \over p}, p}}
\newcommand{\sht}[1]{#1_{{\tiny 1 \over 2}, 2}}
\newcommand{\rad}{r}

\title{A soothing invisible hand:\\ moderation potentials in optimal control}
\author{Debra Lewis\thanks{Mathematics Department, University of California, Santa Cruz, 
Santa Cruz, CA 95064. {\tt lewis@ucsc.edu}. }}

\maketitle

\begin{abstract}
A moderation incentive is a continuously differentiable control-dependent cost term that is identically zero on
the boundary of the admissible control region, and is subtracted from the `do or die' cost function to 
reward sub-maximal control utilization in optimal control systems. A moderation potential is a function on
the cotangent bundle of the state space such that the solutions of Hamilton's equations satisfying appropriate
boundary conditions are solutions of the synthesis problem---the control-parametrized Hamiltonian system 
central to Pontryagin's Maximum Principle. A multi-parameter family of moderation incentives for affinely controlled 
systems with quadratic control constraints possesses simple, readily calculated moderation potentials. One member
of this family is a shifted version of the kinetic energy-style control cost term frequently used in geometric
optimal control. The controls determined by this family approach those
determined by a logarithmic penalty function as one of the parameters approaches zero, while the cost term itself is bounded. 
\end{abstract}


\section{Introduction}
\label{intro}
When modeling a conscious agent, the constant cost function 
of a traditional time minimization problem can be interpreted as representing a uniform stress or 
risk  throughout the task,  while a generalized time minimization cost function models varying
stresses and risks that depend on the current state of the system. Implementation of agent limitations
via cost terms may be more natural---particularly for biological systems---than a possible/impossible dichotomy, 
in which constraints are explicitly incorporated in the state space. For example, consider 
the classic `falling cat' problem, in which a cat is suspended upside down and then released.
(See, e.g.,  \cite{Marey, KS, Mont93}.)  Marey \cite{Marey} gave a qualitative description of the self-righting 
maneuver, supported by high-speed photographs: the cat rotates the front and back halves of its body,
altering the positions of its head and limbs to adjust its moments of inertia, causing the narrowed
half to rotate significantly faster than the thickened half, with zero net angular velocity.
Kane and Scher \cite{KS} introduced a simple mathematical model of a cat,
consisting of a pair of coupled rigid bodies and showed that self-righting with zero angular momentum is possible
without alteration of the moments of inertia of the two halves of the body. To 
rule out the mechanically efficient but fatal solution in which the front and 
back halves simply counter-rotate, resulting in a $360^\circ$ twist in the `cat', Kane and Scher imposed
a no-twist condition in their model. However, actual cats can and do significantly 
twist their bodies, and splay or tuck their limbs; the images in \cite{KS} generated using the mathematical model
and superposed on photographs of an actual cat significantly underestimate the relative motion between the
front and back halves of the body.\footnote{It should be noted that the ultimate goal of Kane and Scher's investigation
was the development of maneuvers that would allow astronauts to alter their orientation without grasping fixed objects
while in zero gravity \cite{life}. The no-twist condition is very reasonable for a spacesuit-clad human---the 
collarbone (not present in cats or many other quadrupeds) limits rotation of the shoulders relative to the hips, facilitating
bipedal motion, and the bulky suit further limits bending and twisting.} Replacing the no-twist condition with a deformation-dependent term in the cost 
function that discourages excessively large relative motions allows more realistic motions.  

The optimal control values for purely state-dependent cost function
lie on the boundary, if any, of the admissible controls set. In some situations, geometric optimization and integration methods can be 
used (see, e.g., \cite{LS, Leimkuhler, synode, LO1}) to solve the restriction of the optimization problem to the boundary of the admissible region. 
If geometric methods are not available or desirable, penalty functions can be used to construct algorithms on an ambient
vector space that respect the boundary due to the prohibitive (possibly infinite) expense of crossing
the boundary; see, e.g., \cite{BG, BB}, and references therein. For some state- and control-dependent cost functions, 
trajectories approaching the boundary of the admissible control region are so extravagant that
the boundary can safely be left out of the mathematical model. However, a close approach to the boundary 
of the admissible region may be appropriate when making the best of a bad situation. Consider again the situation of a cat that
is suspended in mid-air facing upward, then released: the
cat is presumably eager to change its orientation before striking the ground---a typical cat can right 
itself when dropped from heights of approximately one meter. Selection of
an appropriate cost function is essential;
 a very high price for near-maximal control values may yield overly conservative
solutions, while excessively low costs may result in near-crisis responses in almost all situations.
Modeling a system using a family of cost functions parametrized by moderation or 
urgency can reveal qualitative features of optimal solutions that are not readily seen using a
single cost function. 

In time minimization problems, the time required to complete the maneuver is obviously {\it not} specified a priori. 
In more general situations, in which the cost function is a non constant function of state and/or control values,
the duration is still allowed to vary unless specified as fixed.\footnote{Problems in which
the duration is to be determined are given pride of place in Pontryagin et al. \cite{PBGM}:  ``Let us note that (for fixed $x_0$ and $x_1$) the upper
and lower limits, $t_0$ and $t_1$\ldots are not fixed numbers, but depend on the choice of control $u(t)$ which transfers
the phase point from $x_0$ to $x_1$ (these limits are determined by the relations $x(t_0) = x_0$ and $x(t_1) = x_1$).'' 
appears on page 13; the treatment of fixed time problems is deferred to page 66.}
 One of  
Pontryagin's necessary conditions for optimality when the duration is free to vary is the requirement that a trajectory lie within 
the zero level set of the Hamiltonian  \cite{PBGM}. Shifting the Hamiltonian by a constant leaves the 
evolution equations unchanged, but can significantly influence the optimal trajectories via the
initial conditions. Hence formulation of intuitive, geometrically meaningful criteria for the specification
of the constant is central to the analysis of modified time minimization problems.
In \cite{Lewis_cons} we considered generalizations of traditional time minimization problems, allowing both a purely
dependent term  modeling a do-or-die, `whatever it takes' approach, and a
control-dependent term that equals zero on the boundary of the admissible control region; the decreased cost on
the interior can be interpreted as an incentive rewarding sub-maximal control efforts. We introduced two families
of control cost term; one  was modeled on a quadratic control cost (see, e.g.,  \cite{BW},  \cite{Bloch}, \cite{Mont92}, 
and references therein), shifted so as to equal zero on the boundaries
of the admissible regions, while the other, the elliptical moderation incentives, was modeled on the $L_2$ norm
with respect to state-dependent inner products.

Optimal control on nonlinear manifolds has received significant attention in recent years,
particularly situations in which the controls can be modeled as elements of a distribution
within the tangent bundle of the state manifold, corresponding to (partially) controlled 
velocities. See, e.g. \cite{NVdS, Sontaga, Sontagb, Bloch}, and references therein.  Pontryagin's Maximum Principle
on manifolds involves Hamiltonian dynamics with respect to the canonical symplectic structure on the cotangent 
bundle of the state space. We extend the notion of a
moderation incentive to control systems on manifolds and develop general conditions under which moderation incentives
determine a unique optimal control value for each point in the cotangent bundle of the state space. 

One of the advantages of the `kinetic energy' control cost term for systems with boundaryless admissible control regions is that the optimal 
control value is straightforward to compute---it can be computed directly as a non-parametrized Hamiltonian system on $T^* \cS$.
However, if the admissible control regions are bounded, a `kinetic energy' cost term can lead to non-differentiable controls. 
We construct a multi-parameter family of moderation incentives for affine nonlinear control systems with admissible control
regions determined by quadratic forms, and determine the optimal control values for the associated cost functions.
For all but one value of one of the parameters, this family yields differentiable optimal controls on the interior of the admissible control region.
The upper limit of one of the parameter ranges determines generalizations of the shifted kinetic energy' cost, with continuous controls that
fail to be differentiable on the boundary of the admissible control regions. The lower limit, which is not attained, yields controls equal to 
those determined by a traditional logarithmic penalty function. 

Using the optimal control functions, we can construct functions, which we call moderation potentials, on the cotangent bundle of the state space such that the associated Hamiltonian
dynamics are solutions of Pontryagin's synthesis problem. Thus, rather than work with a family of parametrized Hamiltonians, we can find
solutions of a traditional canonical Hamiltonian system. The moderation potentials for the multi-parameter family of moderation incentives
constructed here have a relatively simple form. Some members of the family have particularly simple, geometrically meaningful, expressions.
These members provide a promising alternative to `kinetic energy' cost terms for systems with bounded control regions, allowing use of 
the full range of Hamiltonian machinery.

We illustrate some of the key features of the moderated
control problems using a simple two-dimensional controlled velocity problem: a vertically
launched projectile is guided towards a fixed target; the speed is bounded by a function $\rho$ of the 
horizontal component of the position. The optimal velocity and launch point are to be determined. 
a cost term depending only on the horizontal component 
of the projectile's position models risk from ground-based defense of the target. The cost function is a 
combination of a term depending only on the horizontal component of the projectile's position models risk from ground-based 
defense of the target and a `moderating' function of the control.
We explore the behavior of the solutions of the synthesis problem as the parameters in the cost function are
varied.  

\section{\label{warm_up} Constants matter: moderation incentives}

We first establish notation and context: We assume that the set $\cS$ of possible states is a smooth manifold, and 
consider control problems with state variable $z \in \cS$ and control $u$ in the state-dependent admissible control region $\cA_z$ 
for $z$. We assume that both the state space $\cS$ and the set 
\[
\cA := \lcb (z, u) : z \in \cS \ \mbox{and} \ u \in \cA_z \rcb
\]
of admissible state/control pairs are smooth manifolds. The evolution of the state variable is determined by a controlled vector field $X$. 
Specifically, $\dot z = X(z, u)$ for some continuous map $X: \cA \to T \cS$ satisfying $X(z, u) \in T_z \cS$ for all $(z, u) \in \cA$. 
(Here $T \cS$ denotes the tangent bundle of the state space $\cS$, and $T_z \cS$ denotes the fiber of over $z$ in $T \cS$. See, e.g. \cite{FoM, Bloch}.)
The control problem is to find a duration $\tfin$ and piecewise continuous curve $(z, u): [0, \tfin] \to \cA$ satisfying the 
boundary conditions $z(0)  = z_0$ and $z(\tfin) = z_{\! f}$, with piecewise continuously differentiable state component $z$. 
The optimal control problem with instantaneous cost function $C: \cA \to \R$ is to find the solution that minimizes the total cost 
\[
\int_0^{\tfin} C(z(t), u(t)) dt
\]
over the set of all solutions of the control problem. The fixed time problem is defined analogously, but $\tfin$ is specified. 

Given a purely state-dependent cost function $\hatC: \cS \to \R$, we construct the cost function $C: \cA \to \R$ by subtracting a 
control-dependent term $\tilC(z, u)$ from the unmoderated cost $\hatC(z)$. 
Shifting the cost function by a constant doesn't change the evolution equations of the corresponding Hamiltonian system. However, 
when seeking solutions of the optimization problem, shifting the Hamiltonian by a constant can have a significant effect
via the condition that the Hamiltonian equal zero along solutions of the synthesis problem. To guide the selection of the control-dependent
function $\tilC$, we regard that term as an incentive for sub-maximal control investment, rather a penalty. This motivates the condition that 
$\tilC(z, u) = 0$ if $u \in \partial \cA_z$.  

\begin{Definition}
Given an admissible space $\cA$, $\tilC \in \cC^1(\cA, [0, \infty))$  is a {\em moderation incentive} for $\cA$ if  for all $z \in \cS$
$u \in \partial \cA_z$ implies $\tilC(z, u) = 0$. 

If there are continuous functions $q: \cA \to [0, 1]$ and $\Phi: \cS \times [0, 1] \to \R$ with $\Phi^{-1}(0) = \cS \times \{1 \}$
such that for every $z \in \cS$,
$\partial \cA_z = \{ u \in \cA_z: q(z, u) = 1 \}$
and $s \mapsto \Phi(z, s)$ is a decreasing function,
then $\tilC(z, u) := \Phi(z, q(z, u))$ is a {\em monotonic moderation incentive} for $\cA$ and $q$. 
\end{Definition}

The following example illustrates the influence of shifting the cost function of an optimal control system for which the time is not fixed. 
We consider a two dimensional system with controlled velocities. Starting
from the horizontal axis, with vertical initial velocity, the goal is to hit a 
target $(\xfin, \yfin)$. We assume that the projectile starts to the right of the target and
consider a non-increasing unmoderated position-dependent cost term $\hatC: [\xfin, \infty) \to \R^+$ 
depending only on the horizontal component of the position, modeling risk due to ground-based 
defense of the target, combined with a control- (and possibly position-) dependent moderation term.
Given the final height $\yfin$, we seek smooth  trajectories $(x, y): [0, \tfin] \to \R^2$ satisfying 
$y(0) = \dot x(0) = 0$, $x(\tfin) = \xfin$, and $y(\tfin) = \yfin$. Neither the
launch point $(x_0, 0)$ nor final time $\tfin$ are specified in advance. We have direct control over
the velocity, with the constraint that the speed of the projectile never exceeds one.

We consider 
a pair of two-parameter families of cost functions, differing only by a constant. One parameter, $c$, scales a purely position-dependent term; 
the second, $\mu$, scales a control-dependent term. One family yields inflexible solutions---the solution path is entirely 
determined by the boundary conditions, while the speed is simply rescaled by the ratio of the two parameters. The other family, in which the
parameter $\mu$ scales a moderation incentive, has solutions for
which the optimal path and speed both depend nontrivially on the parameters $c$ and $\mu$.  Here we summarize some of the key
features of this system---our intent is only to remind the reader that analogous choices can have a profound 
influence on the solutions of optimal control systems if the total time is not specified {\it a priori}, and hence should be systematically selected. 
A generalization of this system is analyzed in \S \ref{vto}.

The cost functions
\beq{vto_costs1}
C_{\rm ke}(x, y, \dot x, \dot y; \mu) := \frac \mu 2  \norm{(\dot x, \dot y)}^2 +  \frac c {2 \, x^2} 
\eeq
and
\beq{vto_costs2}
C_{\rm mi}(x, y, \dot x, \dot y; \mu) := C_{\rm ke}(x, y, \dot x, \dot y; \mu)  + 1 - \frac \mu 2.
\eeq
differ only by a constant, but there are important differences in the behavior of the solutions. 
The Hamiltonians associated to $C_{\rm ke}$ and $C_{\rm mi}$ (see \S \ref{modpot}) equal that of a 
point mass with mass $\mu$ and potential energy $- \frac c {2 \, x^2} + \mbox{constant}$. The solutions trace
 out segment of ellipses with principal axes $\sqrt{1 - (\xfin/x_0)^2}$ and $\yfin/x_0$. The constant influences 
the  solutions of the synthesis problem through the condition that the Hamiltonian equal zero on an optimal trajectory. 

\begin{figure}[t]
\begin{center}
\includegraphics[width=1.6in]{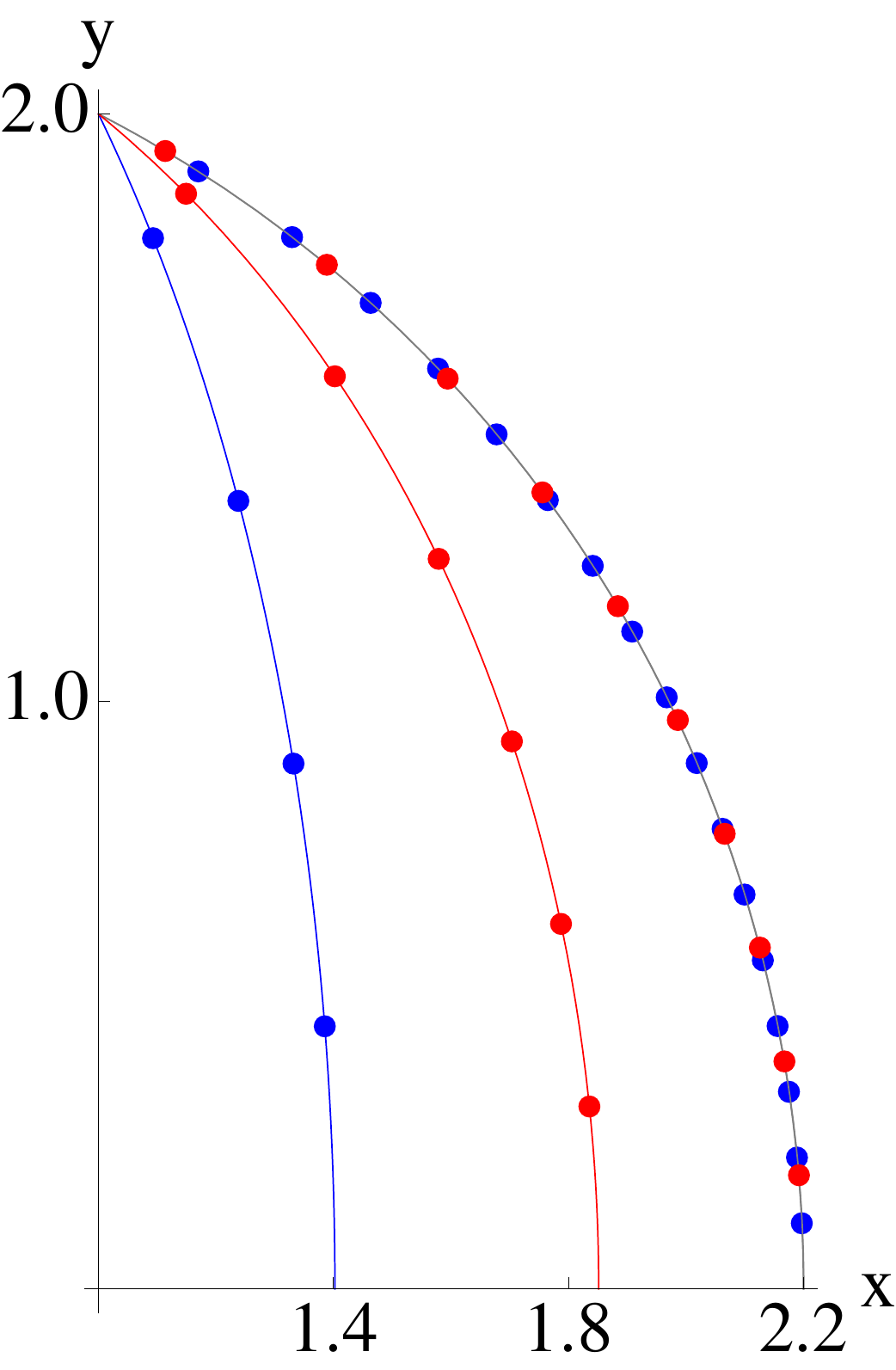} \qquad 
\includegraphics[width=1.6in]{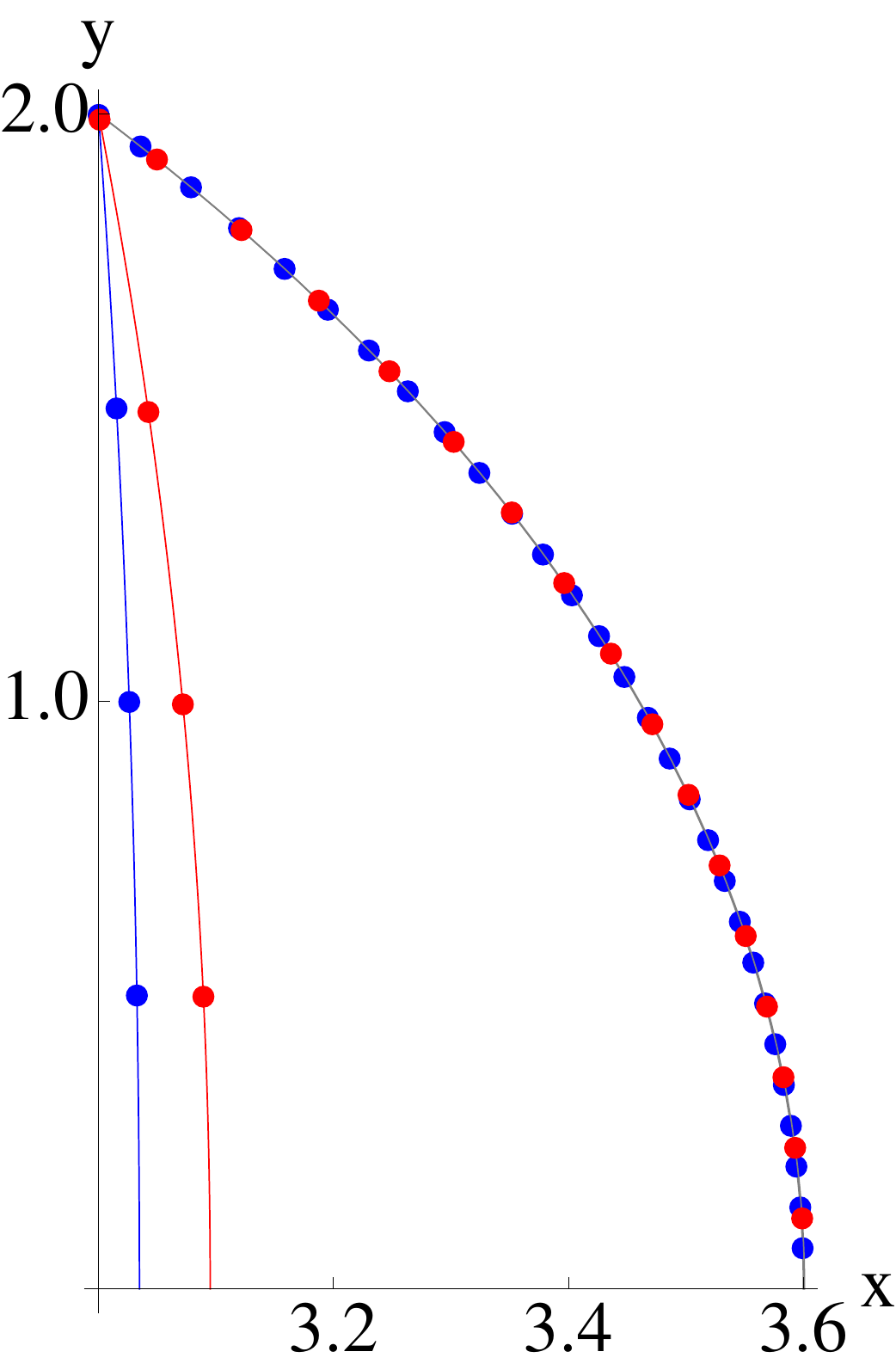} 
\end{center}
\caption{\label{vto_plots}Solutions of the vertical take-off targeting problem problem for sample target positions $(\xfin, 2)$,
cost function $C_{\rm mi}$, and defense strengths $c$. 
Blue: $c = \frac 1 2$, red: $c = \frac 3 2$; left: $\xfin = 1$, right: $\xfin = 3$.
Dots indicate projectile position at times $t = \frac j 2$, $j \in {\Bbb N}$; colored lines
indicate the traces of solutions for $\mu = \mu_{\rm min}$; solid grey lines indicate those for $\mu = \mu_{\rm max} = 2$. 
All solutions for the cost function $C_{\rm ke}$ lie on the grey curves.
}
\end{figure}

\begin{itemize}
\item
An optimal solution for $C_{\rm ke}$, traces an arc of a circle centered at the origin.
The parameters $c$ and $\mu$ influence the solution only through the rescaling of the speed by $\sqrt{c/\mu}$. 
\item
The optimal starting point $x_0$ for a solution minimizing $C_{\rm mi}$
depends nontrivially on both $c$ and $\mu$; specifically,
\[
x_0^2 = \smallfrac b 2 + \sqrt{\lp \smallfrac b 2 \rp^2 + a \, d}
\qquad \mbox{where} \qquad a = \xfin^2 + \yfin^2, \quad d = \smallfrac c {2 - \mu}, \quad b = \xfin^2 - d.
\]
$x_0$ and $\tfin$ are increasing functions of  $\mu$. Since smooth solutions satisfying the
control constraint $1 \geq \dot x^2 + \dot y^2$ exist only if 
$1 + \frac c {2 \, \xfin^2} = \mu_{\rm min} \leq \mu_{\rm max} = 2$,
we must have $\xfin^2 \geq \frac c 2$.  The projectile follows a circular arc when $\mu = \mu_{\rm max} = 2$.
\end{itemize}

\section{Affine nonlinear control systems and ellipsoidal admissible control regions}

In \cite{Lewis_cons} we introduced the notion of a moderation incentive for control systems in which both the state space and control regions were subsets 
of $\R^n$ and $\R^k$, and focused on the situation $n = k$. Some of the incentives considered yielded simple explicit expressions for the optimal control values as rescalings of the auxiliary variable (the Lagrange multiplier in the Pontryagin formulation of the control problem). We now extend that 
strategy to nonlinear manifolds and introduce a family of moderation incentives for systems in which the admissible control
regions $\cA_z$ are the unit balls with respect to state-dependent norms. The optimal controls for these incentives are rescalings of the
image of the auxiliary variable (now an element of the cotangent bundle of the state manifold) under a mapping determined by the norms.

\begin{Definition}
Given a family of positive-definite quadratic forms $Q_z$ on $\R^k$, such that $z \mapsto Q_z$ is $\cC^1$, we will say 
that a control problem with admissible region 
\beq{elliptic_control_region}
\cA := \lcb (z, u)  \in \cS \times R^k :  Q_z(u) \leq 1 \rcb
\eeq
has {\rm ellipsoidal control regions}.
\end{Definition}

If there are continuous vector fields $f$ and $g_j$, $j = 0, \ldots, k$,  on $\cS$ such that 
\[
\dot z = X(z, u) = f(z) + \Sigma_{j = 1}^k u_j g_j(z), 
\]
the system is said to be {\it affine nonlinearly controlled}; $f$ is the {\it drift} vector field. (See, e.g., \cite{Bloch, Sontagb}.)

Given a control system with ellipsoidal control regions and affine controls, for each $z \in \cS$, 
let $L_z$ and $\dla \ , \ \dra_z$ denote respectively the invertible symmetric linear map and inner product on 
$\R^k$ satisfying 
\[
Q_z(u) = \la u, L_z^{-1} u \ra
\sands
\dla u, v \dra_z =  \la u, L_z^{-1} v \ra
\]
for all $u \in \R^k$,
and define the maps $M_z: \R^k \to T \cS$, $\lambda: T^* \cS \to \R^k$, and $\ell: T^* \cS \to [0, \infty)$ by
\beq{lambda_def}
M_z u := \Sigma_{j = 1}^k u_j g_j(z), \qquad \qquad
\lambda(\ps) := L_z (M_z^* \ps),
\eeq
and
\beq{tilq_def}
\ell(\ps)^2 := Q_z(\lambda(\ps)) = \ps \cdot (X(z, \lambda(\ps)) - X_0(z)).
\eeq
(Here $T^* \cS$ denotes the cotangent bundle of $\cS$; see, e.g., \cite{Bloch, FoM}.)
Finally, define the map $\nu: \ell^{-1}(\R^+) \to \R^k$ by 
\[
\nu(\ps) := \smallfrac 1 {\ell(\ps)} \lambda(\ps) \in \partial \cA_z.
\]

\begin{Proposition}
\label{ellip_control}
Consider a control system with ellipsoidal control regions and affinely controlled evolution. 
Let $F \in \cC^0(\cS \times [0, 1], [0, \infty))$ 
be a function satisfying $F^{-1}(0) = \cS \times \{0\}$ and such that
for every $z \in \cS$, $F_z(x) := F(z, x)$ is increasing and differentiable on $(0, 1)$,
with $\lim_{x \to 1}F_z'(x) < \infty$.
Given $p \geq 1$, define $x: \cS \times \R^+ \times [0, 1]  \to \R$ by
 \beq{x_def}
x(z, \ell, s) :=  \ell \, s + F\! \lp z, 1 - s^p \rp. 
\eeq
If there is a function $\sigma$ on $\cS \times \R^+$ such that for every $(z, \ell) \in \cS \times \R^+$
$s \mapsto x(z, \ell, s)$ achieves its maximum exactly at $\sigma(z, \ell)$, then 
the moderation incentive
\beq{gen_c_def}
\tilC(z, u) := F \lp z, 1 - Q_z(u)^{\frac p 2} \rp
\eeq
has optimal control value 
\beq{gen_opt_u}
\upsilon(\ps) := \lcb \begin{array}{ll}\sigma(\ps) \,\nu(\ps) \qquad & \ell(\ps) \neq 0 \smallskip \\
0  & \ell(\ps) = 0 \end{array} \right .
\eeq
at $\ps \in T^* \cS$.  
\end{Proposition}

\proof
Fix $ \ps \in T^* \cS$. Define $F_z: [0, 1] \to [0, \infty)$ by $F_z(x) := F(z, x)$, $\tilfps: \cA_z \to \R$ by
\[
\tilfps(u) :=  \ps \cdot (X(z, u) - f(z)) + F_z  \lp 1 - Q_z(u)^{\frac p 2} \rp = \dla \lambda(\ps), u \dra_z + F_z  \lp 1 - Q_z(u)^{\frac p 2} \rp,
\]
and $x_{\ps}: [0,1] \to \R$ by $x_{\ps}(s) := x(z, \ell(\ps), s)$.

If $\lambda(\ps) = 0$, $\tilfps(u) = F_z  \lp 1 - Q_z(u)^{\frac p 2} \rp$, which achieves its maximum at $u = 0$.

We now show that if $\lambda(\ps) \neq 0$, then $\tilfps$
takes its maximum on the line segment
\[
\{ s \, \nu(\ps) : 0 \leq s \leq 1 \},
\]
and hence, since $\tilfps(s \, \nu(\ps)) =  x_{\ps}(s)$,  the maximum of  $\tilfps$ 
coincides with the maximum of $x_{\ps}$.
The restriction of  $\tilfps$ to the interior of $\cA$ is differentiable, with gradient 
\[
\nabla \tilfps(u) = \lambda(\ps) - p \, F_z' \!\lp 1 - Q_z(u)^{\frac p 2} \rp Q_z(u)^{\frac p 2 - 1}.
\]
Hence if $\lambda(\ps) \neq 0$, any critical point of $\tilfps$ in the interior of $\cA$ has the form 
$u = s  \, \nu(\ps)$ for $s \in (0, 1)$ satisfying 
\beq{s_eq}
\frac {\ell(\ps)} p  = s^{p - 1} \, F_z'\!\lp 1 - Q_z(s  \, \nu(\ps))^{\frac p 2} \rp = s^{p - 1} \, F_z'\! \lp 1 - s^p \rp.  
\eeq
Note that if $\lambda(\ps) \neq 0$,  and hence $\ell(\ps) \neq 0$, then $s$ satisfies 
\eqr{s_eq} iff $s$ is a critical point of $x_{\ps}$. 

Since $F(z, 0) = 0$ implies that $F_z  \lp 1 - Q_z(u)^{\frac p 2} \rp = 0$ for $u \in \partial \cA_z$, 
\[
\max_{u \in \partial \cA_z} \tilfps(u) 
= \max_{Q_z(u) = 1}  \dla \lambda(\ps), u \dra_z
\]
Hence a standard Lagrange multiplier argument shows that the restriction of  $\tilfps$ to $ \partial \cA_z$
achieves its maximum, $\ell(\ps) = x_{\ps}(1)$, at $\nu(\ps)$. 

Finally, $\tilfps(0) = F_z(1) = x_{\ps}(0)$. 
\prfend

\rmk
\label{scale_mu}
If \eqr{gen_opt_u} is the optimal control for a moderation incentive of the form \eqr{gen_c_def}, with scaling factor $\sigma$, then
given $\mu \in \cC^0(\cS, \R^+)$, 
\[
\tilC_\mu(z, u) := \mu(z) \, F \lp z, 1 - Q_z(u)^{\frac p 2} \rp
\]
is a moderation incentive with scaling factor obtained by replacing $\ell(\ps)$ with 
$\ell_\mu(\ps) := \frac {\ell(\ps)}{\mu(z)}$
in \eqr{gen_opt_u}.
\prfend

A moderation incentive is required to take the value zero on the boundary of the admissible control regions, but is not required to have a finite derivative
there. If the limit of the derivative of the incentive  as $\partial \cA_z$ is approached is infinite, the optimal control must lie in the interior of $\cA_z$.  

\begin{Lemma}
\label{easy_case}
If $F: \cS \times [0, 1] \to [0, \infty)$ satisfies $F^{-1}(0) = \cS \times \{0\}$, and for every $z \in \cS$ the function $F_z(x) := F(z, x)$ is $\cC^2$ on $(0, 1)$, with decreasing positive derivative satisfying 
\beq{limit_derivs}
 \lim_{x \to 0}F_z'(x) = \infty,
\eeq
then $s \mapsto x(z, \ell, s)$ given by (\ref{x_def}) achieves its maximum at a unique point 
$s_*(z, \ell) \in \lp 0, 1 \rp$ if $p > 1$, or if $p = 1$ and $F_z$ is strictly decreasing.

The associated map $\sigma: \ell^{-1}(\R^+) \to (0, 1)$ given by $\sigma(\ps) := s_*(z, \ell(\ps))$ is $\cC^1$.
\end{Lemma}
\proof
Setting $y = 1 - s^p$ and $c = \frac {\ell(\ps)} {p} > 0$, \eqr{s_eq} takes the form 
$c (1 - y)^{\frac 1 p - 1} = F_z'(y)$, with unique solution $y(c) \in (0, 1)$.
The map $s \mapsto s^{p - 1} \, F_z'\! \lp 1 - s^p \rp$ has a $\cC^1$ strictly positive derivative on $(0, 1)$.
Hence the Implicit Function Theorem implies that the map $\sigma$ determined by \eqr{s_eq} is $\cC^1$
on $\ell^{-1}(\R^+)$. 
\prfend

We now define a family of monotonic moderation incentives for affinely controlled systems with quadratic control costs. These 
generalize the moderation incentives introduced in \cite{Lewis_cons}. The `dogleg parameters' $\alpha \in (0, 1]$ and $p \geq 1$
can be interpreted as tuning the overall shape of the control response curve, while the state-dependent moderation strength function 
$\mu \in \cC^0(\cS, \R^+)$ scales the instantaneous control cost. (Use of the term `dogleg' is motivated by the shape of the response
curve for values of $\alpha$ near 1; varying these parameters alters the abruptness of the dogleg bend.)

\begin{figure}[t]
\begin{center}
\includegraphics[width=3in]{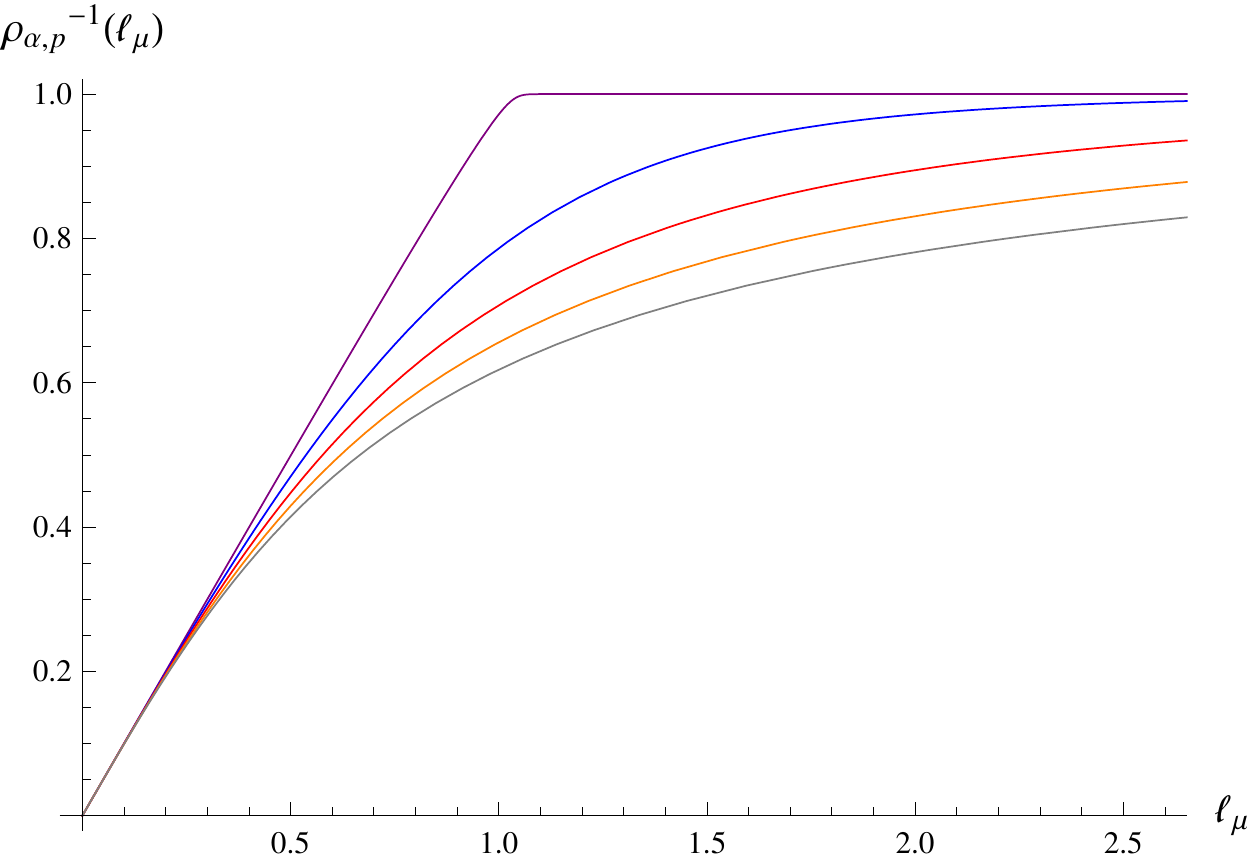}
\includegraphics[width=3in]{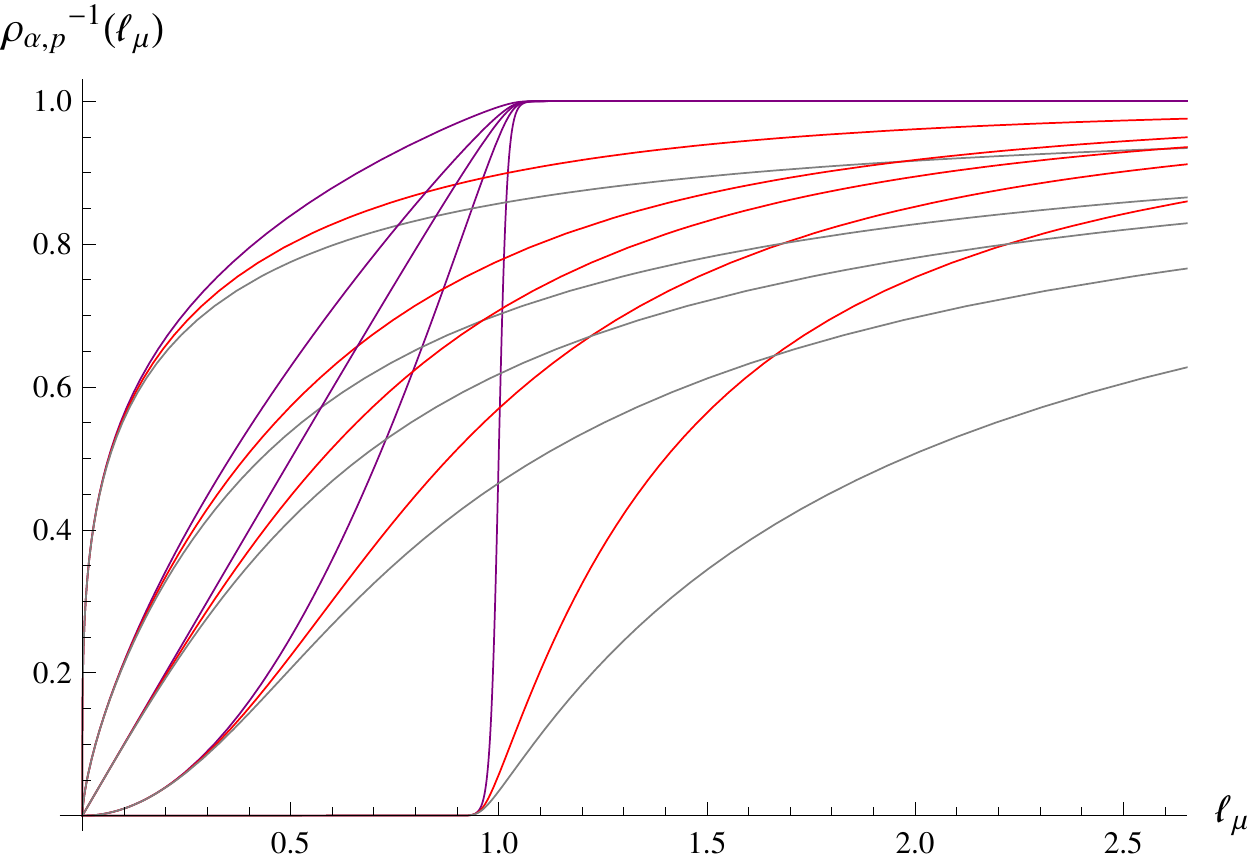}
\end{center}
\caption{\label{opti_plots}
Plots of $\rho_{\alpha, p}^{-1} \lp \ell_\mu \rp$ for different values of the dogleg parameters $\alpha$ and $p$. 
Purple: $\alpha = \frac {99}{100}$; blue: $\alpha = \frac 3 4$; red: $\alpha = \half$; orange:
$\alpha = \frac 1 4$; gray: limiting case $\alpha \to 0$. 
Left: $p = 2$. Right: $p = 1.01, 1.5, 2, 2.5, 5$; convexity for small values of $s$ increases with $p$. 
(The approximate step function in the right hand graph
 is associated to $\alpha = \frac {99}{100}$, $p = 1.01$.)
}
\end{figure}

\begin{Theorem}
\label{emi}
Given $0 < \alpha \leq 1\leq p$, excluding $\alpha = 1 = p$, and $\mu \in \cC^0(\cS, \R^+)$,
\beq{ellip_mod_incentive}
\tilC_{\alpha, p}(z, u; \mu) = \smallfrac {\mu(z)} {p \, \alpha}  \lp 1 - Q_z(u)^{\frac p 2} \rp^\alpha
\eeq
is a monotonic moderation incentive for $\cA$. 

If $0 < \alpha < 1$, the unique optimal control parameter  \eqr{gen_opt_u} for $\tilC_{\alpha, p}$ 
has the scaling
\beq{sigma_alpha}
\sigma_{\alpha, p}(\ps; \mu) = \rho_{\alpha, p}^{-1}(\ell_\mu(\ps)),
\eeq
where $\rho_{\alpha, p}: [0, 1] \to [0, \infty)$ and $\ell_\mu: T^* \cS \to [0, \infty)$ are given by
\[
\rho_{\alpha, p}(s) := s^{p - 1} \lp 1 - s^p \rp^{\alpha - 1}
\sands \ell_\mu(\ps) := \frac {\ell(\ps)}{\mu(z)} .
\]
If $\alpha = 1 < p$, then 
\beq{QCC_ocp}
\sigma_{1, p}(\ps; \mu) : =  \min \lcb \ell_\mu(z)^{\frac 1 {p - 1}} , 1 \rcb
\eeq
is the optimal scaling.
\end{Theorem}
\begin{proof}
For $0 < \alpha \leq 1$, $F_\alpha(x) := \frac 1{\alpha} x^\alpha$ is differentiable, with decreasing positive derivative, on $(0, 1]$. 
For $0 < \alpha < 1$,  $\lim_{x \to 0}F_\alpha'(x) = \infty$, so the rescaling of $F_\alpha$
by $\frac {\mu(z)} p$ satisfies the conditions of Lemma~\ref{easy_case}.

The case $\alpha = 1< p$ requires a direct application of Proposition \ref{ellip_control}, since $F_1' \equiv  1$.
In this case, 
\[
\frac {x_{\ps}(s)}{ \mu(z)} =  \ell_\mu(\ps) \, s + \smallfrac 1 p \lp 1 - s^p \rp 
\]
is the restriction of a polynomial to $\lsb 0, 1 \rsb$. If $\ell_\mu(\ps) \leq 1$, the maximum of 
$x_{\ps}$ coincides with that of the polynomial, which occurs at $s = \ell_\mu(\ps)$.  
If $\ell_\mu(\ps) \geq 1$, the maximum occurs at one of the endpoints;  since 
\[
x_{\ps}(0) = \frac {\mu(z)} p < \ell(\ps) = x_{\ps}(1),
\]
in this case $x_{\ps}$ achieves its maximum of $\ell(\ps)$ at $1$.
\end{proof}

\begin{figure}[t]
\begin{center}
\includegraphics[width=3in]{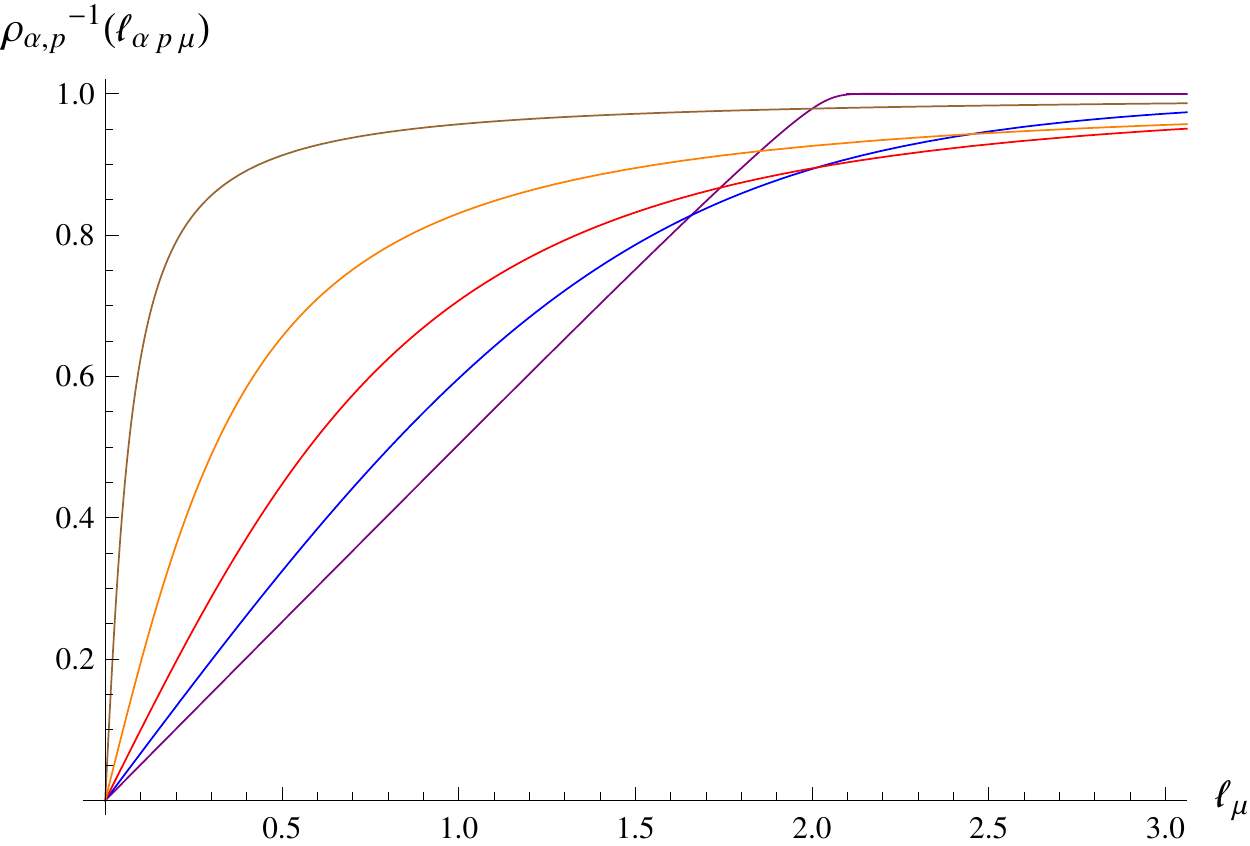}
\includegraphics[width=3in]{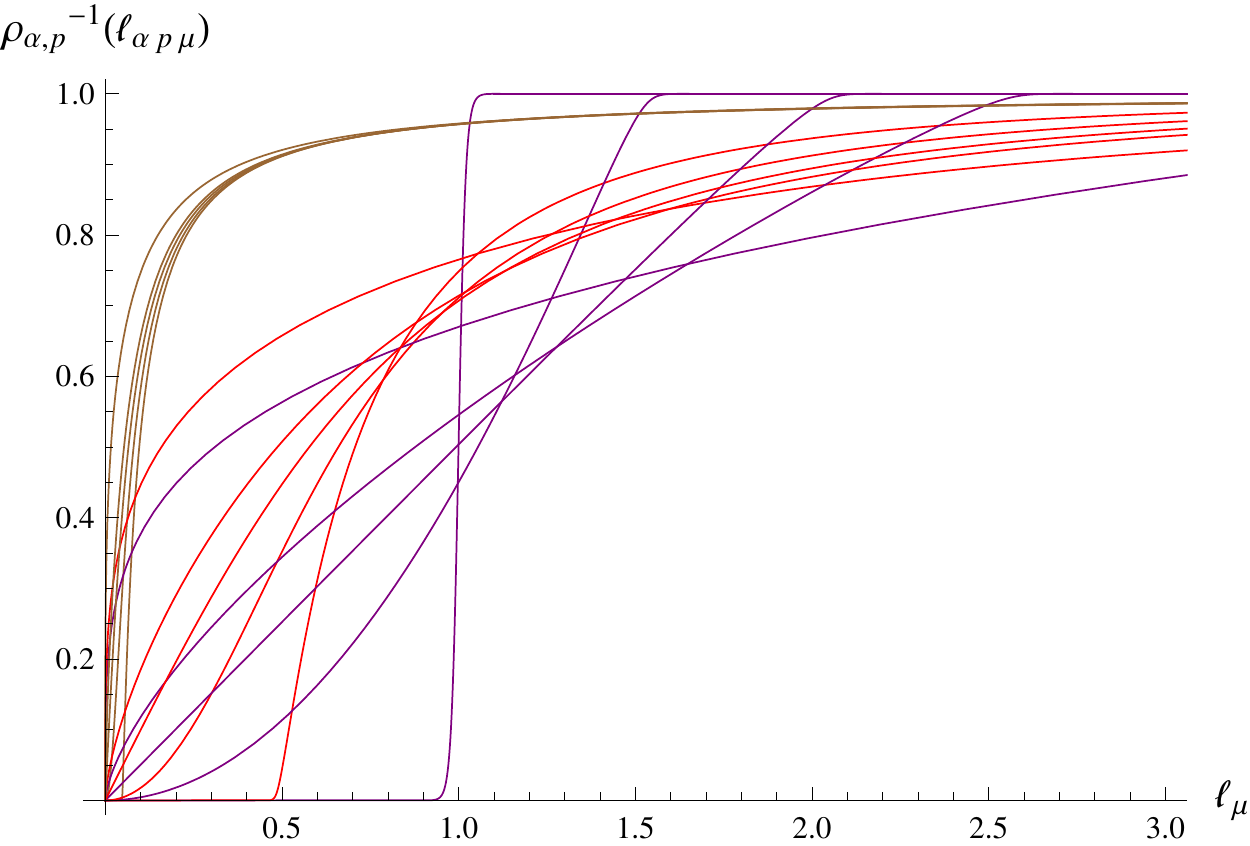}
\end{center}
\caption{\label{opti_plots_scaled}
Plots of $\rho_{\alpha, p}^{-1} \lp \frac {\ell_\mu}{\alpha \, p} \rp$ for different values of the dogleg parameters $\alpha$ and $p$. 
Purple: $\alpha = \frac {99}{100}$; blue: $\alpha = \frac 3 4$; red: $\alpha = \half$; orange:
$\alpha = \frac 1 4$; brown: $\alpha = \frac 1 {20}$. 
Left: $p = 2$. Right: $p = 1.01, 1.5, 2, 2.5, 5$.}
\end{figure}

\subsection{Special cases: $\alpha = 1$, $\alpha \to 0$, and $\frac 1 \alpha = p$ }

For some special values of the parameters $\alpha$ and $p$, simple closed form expressions for the 
optimal scaling $\sigma_{\alpha, p}$ exist. 

The dogleg parameter values $\alpha = 1$, $p = 2$ correspond to the widely used `kinetic energy' style control cost, which is 
used in the targeted attack problem in \S \ref{intro} and \S \ref{vto}. Note that when $\alpha = 1$, the optimal scaling 
is not differentiable at $\partial \cA_z$.

The limit $\lim_{\alpha \to 0} \tilC_{\alpha, p}$ is not well-defined, but 
\[
\rho_{0, p}(s) := \lim_{\alpha \to 0} \rho_{\alpha, p}(s) = \frac {s^{p - 1}} {1 - s^p} = - \smallfrac {d \ }{ds}  \ln \lp \lp1 - s^p \rp^{\frac 1 p} \rp
\]
is well-defined and invertible on $[0, 1)$. In particular, the optimal scaling associated to the logarithmic control cost 
\[
\breve C(z, u) := - \half \, \ln (1 - Q_z(u))
\]
is $\rho_{0, 2}^{-1}(\ell_\mu(\ps))$. Thus the controls determined by the family $\tilC_{\alpha, 2}$ determine a homotopy 
between optimal controls determined by a `kinetic energy'  control cost and a logarithmic `penalty function' cost.
(Logarithmic penalty functions are widely used in the engineering literature to enforce inequality constraints.)

In the case $\alpha = \frac 1 p < 1$, we can explicitly invert $\shp{\rho}$.
\begin{Cor}
\label{qpone}
\beq{qp_case}
\shp{\sigma}(\ps; \mu) = \lp 1 + \ell_\mu(\ps)^{-q} \rp^{- \frac 1 p} \qquad \mbox{for} \qquad  \frac 1 p + \frac 1 q = 1,
\eeq
and hence
\beq{qp_opt}
\shp{\upsilon}(\ps; \mu) =  \frac {\ell(\ps)^{q - 2}}{\lp \mu(z)^q + \ell(\ps)^q \rp^{\frac 1 p}}  \lambda(\ps)
\eeq
if $\lambda(\ps) \neq 0$.
\end{Cor}

\begin{figure}[t]
\label{opt_proj_hemi}
\begin{center}
\includegraphics[width=3.75in]{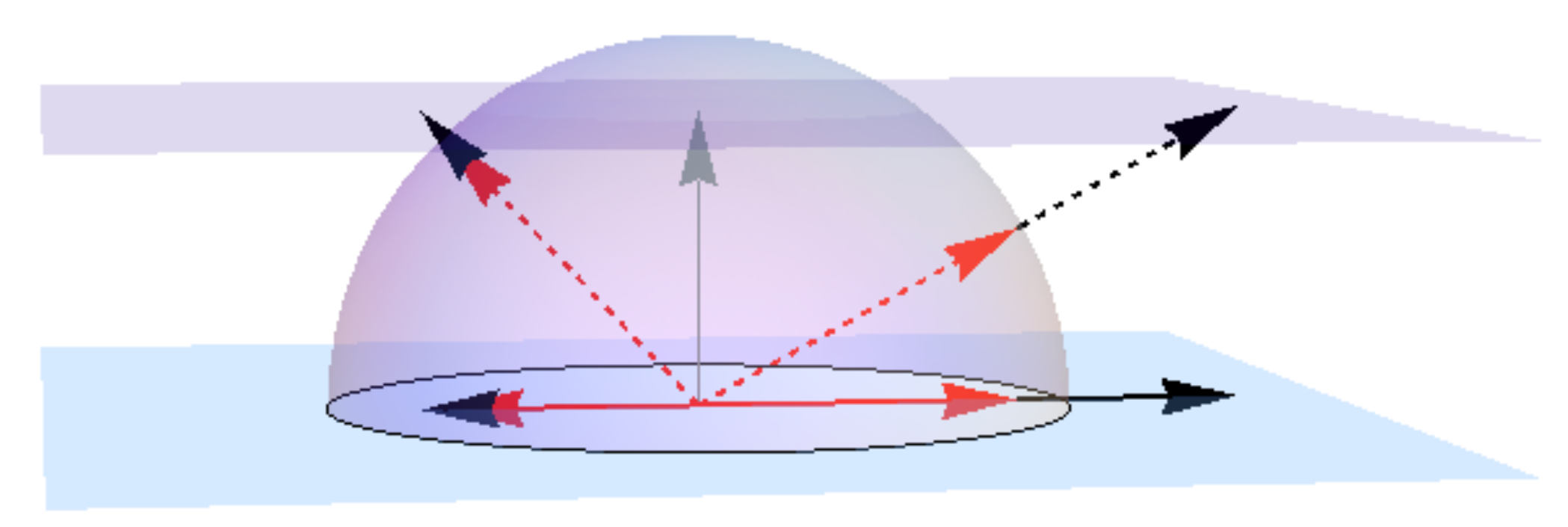}
\end{center}
\caption{\label{opti_proj_hemi}
The optimal control for $\frac 1 \alpha = p = 2$ implemented as lift into control--moderation space, followed by projection onto the unit sphere, then
projection back into control space. Black solid arrows: $\lambda(\ps)$; black dashed: $(\lambda(\ps), \mu(z))$; red dashed: 
$\frac {(\lambda(\ps) , \mu(z))}{\norm{(\lambda(\ps) , \mu(z))}_{Q_z}}$; red solid: $\frac {\lambda(\ps)}{\norm{(\lambda(\ps) , \mu(z))}_{Q_z}}$.
}
\end{figure}

\rmk
When the drift field is trivial,  the optimal rescaling has the following geometric interpretation: the optimal control $\sht{\upsilon}$ satisfies
\[
\upsilon_{\half, 2}(\ps)
= \frac {\lambda(\ps)} {\norm{(\lambda(\ps) , \mu(z))}_{Q_z}},  
\]
where $\norm{(u, t)}_Q^2 = Q(u) + t^2$ is the norm on $\R^{k + 1}$ induced by a quadratic form $Q$ on $\R^k$. Thus
$\upsilon_{\half, 2}(\ps)$ is the control component of the projection of $(\lambda(\ps) , \mu(z))$ onto the $\norm{\ }_Q$ unit ball in $\R^{k + 1}$.
(See Figure \ref{opt_proj_hemi}.) We will further investigate the moderated controls for $\alpha \, p = 1$, particularly that for $p = 2$, in future work. 
\rmkend

\section{Moderation potentials and the synthesis problem}
\label{modpot}

Pontryagin's Maximum Principle relates optimal control to Hamiltonian dynamics:  
If the state space $\cS$ is an $n$-dimensional subset of $\R^n$, then associated to a solution
 $(z, u): [0, \tfin] \to  \cA \subseteq \R^n \times \R^k$ of the control problem  
minimizing the total cost there is a curve $\lambda: [0, \tfin] \to \R^n$ such that the curve $(z, \lambda)$
satisfies Hamilton's equations for the time-dependent Hamiltonian
\[
H_t(\tilde z, \lambda) := \lambda^T X(\tilde z, u(t)) - C(\tilde z, u(t)),,
\]
and
\[
H_t(z(t), \lambda(t)) = \max_{u \in \cA_{z(t)}} \lp \lambda(t)^TX(z(t), u) - C(z(t), u) \rp.
\]
(See \cite{PBGM} for the precise statement and proof of the Maximum Principle.) 
Pontryagin's optimality conditions are necessary, but not sufficient. Their appeal lies in their
constructive nature---well-known results and techniques for boundary value problems and Hamiltonian 
dynamics can be used to construct the pool of possibly optimal trajectories. 
This construction is referred to as the synthesis problem in \cite{PBGM}; we will make use of that terminology here.) 

The generalization of Hamilton's equations to a nonlinear state manifold $\cS$ utilizes the canonical symplectic structure on the 
cotangent bundle $T^* \cS$ of the state manifold. See, e.g., \cite{Bloch} for additional background and discussion. We now
introduce the formulation of the synthesis problem that will be used here. Given our focus on systems with state-dependent
admissible control regions, we introduce the condition that the parametrized Hamiltonians be extendable to neighborhoods
of possibly optimal values. 

Given that the admissible control regions can vary with the state variable, we explicitly require that the vector field $X$ and cost term $C$
with fixed control value be extensible to neighborhoods of the points of interest. Let $\pi: T^* \cS \to \cS$ denote the canonical projection, with 
$\pi^{-1}(z) = T_z^* \cS$, 
\[
\cP := \{ (\ps, u) : (\pi(\ps), u) \in \cA \},
\]
and $\bbP_1: \cA \to \cS$ denote projection onto the first factor. Given $F \in \cC^1(T^* \cS)$, let $X_F$ denote the Hamiltonian vector
field determined by $F$ and the canonical symplectic structure on $T^* \cS$. Define $H: \cP \to \R$ and $\chi: T^* \cS \to \R$ by
\beq{chi_def}
H(\ps, u) := \ps \cdot X(z, u) - C(z, u) 
\sands \chi(\ps)  := \max_{u \in \cA_z} H(\ps, u).
\eeq

\begin{Definition}
If for every $(\ps, u_*) \in \cP$  satisfying $H(\ps, u_*) = \chi(\ps)$
there is a neighborhood $\cV$ of $z$ such 
that the restrictions of $X( \ \cdot \ , u_*)$ and $C( \ \cdot \ , u_*)$ to $\cV$ are restrictions to $\cV \cap \bbP_1(\cA) $ of $\cC^1$ 
maps on $\cV$, then $X$ and $C$ are {\em synthesizable}.

A curve $(\Psi, \upsilon): [0, \tfin] \to \cP$ satisfying 
\[
\dot \Psi(t) = X_{H_t}(\Psi(t))
\sands
 H_t(\Psi(t)) = \chi(\Psi(t))
 \]
for the local time-dependent Hamiltonians $H_t(\ps) :=  H(\ps, \upsilon(t))$ determined by synthesizable $X$ and $C$
 is a solution of the {\rm synthesis problem} determined by $X$, $C$ and the boundary data $z_0 = \pi(\Psi(0))$ and $\zfin = \pi(\Psi(\tfin))$
 if $H_0(\Psi(0)) = 0$.
 
If $(\Psi, \upsilon)$ satisfies the all of the above conditions except the condition that $H_0(\Psi(0)) = 0$, then $(\Psi, \upsilon)$ is a  solution of the 
{\rm fixed time synthesis problem} of duration $\tfin$. 
\end{Definition}

We will focus on finding solutions of the synthesis problem, and will not formulate general conditions under which such solutions
are in fact global minimizers. 

One of the advantages of the `kinetic energy' control cost term widely used in geometric optimal control
for systems with controlled velocities and unbounded admissible control regions is that the optimal 
control value is straightforward to compute---it is simply the inverse Legendre transform of $\ps$---and hence
solutions of the synthesis problem can be computed directly as a non-parametrized Hamiltonian system on $T^* \cS$.
However, if the admissible control regions are bounded, a `kinetic energy' cost term can lead to non-differentiable controls. (See, e.g. 
\cite{Lewis_cons}.) We now identify conditions under which a control-dependent cost function determines solutions 
of the synthesis problem corresponding to solutions of a traditional Hamiltonian system

Nondegeneracy of the symplectic structure guarantees that two Hamiltonian vector fields agree at $\ps$ iff $\ps$ is a critical point
of the difference of the Hamiltonians. For fixed $\ps \in T^* \cS$, we can define $h_{\ps} \in \cC^1(\cA_z, \R)$ by $h_{\ps}(u) := H(\ps, u)$.
If $H$ is $\cC^1$ and $h_{\ps}$ achieves its maximum at
$u_* \in \cA_z^o$, then $u_*$ is a critical point of $h_{\ps}$. It follows that if $H$ is $\cC^1$ and there is a $\cC^1$ map $\upsilon$ such that 
$H(\ \cdot \ , \upsilon(\ \cdot \ )) = \chi$ and $\upsilon(\ps) \in \cA_z^o$ for every $\ps \in T^* \cS$, then 
\[
{\rm d}(\chi - H(\ \cdot \ , u_*))(\ps)(w_{\ps}) = d H(\ps, u_*)(0,d_{\ps} \upsilon(w_{\ps})) 
= ({\rm d}h_{\ps}(u_*)(d_{\ps} \upsilon(w_{\ps}))) = 0
\]
for $u_* = \upsilon(\ps)$ and all $w_{\ps} \in T_{\ps} T^* \cS$. 

If there is a $\cC^0$ map $\upsilon$ such that $H(\ \cdot \ , \upsilon(\ \cdot \ )) = \chi$, but $\upsilon$ is not everywhere differentiable, or
can take values on the boundaries of the admissible control region, then the above argument  is not applicable, but we may still be able
to replace the parametrized Hamiltonians $H(\ \cdot \ , u_*)$ with the function $\chi$.

The solution of the synthesis problem can be simplified given a feedback law that allows replacement of the control-parametrized Hamiltonian
with a conventional autonomous Hamiltonian on the cotangent bundle $T^* \cS$. We now show that the moderation incentives $\tilC_{\alpha, p}$
have such control laws. The key concerns are 
formulation of relatively simple expressions for the Hamiltonians and verification of the differentiability of the Hamiltonian on the boundaries
of the admissible control regions.
 
 \begin{Proposition}
 \label{build_sols}
If  
\begin{enumerate}
\item
$X$ and $C$ are synthesizable, 
\item
$\chi$ given by \eqr{chi_def} is $\cC^1$,  
\item
$\Psi: [0, \tfin] \to T^* \cS$ is a solution of the canonical Hamiltonian system with Hamiltonian $\chi$, and 
\item
there is a curve $u: [0, \tfin] \to \cP$ such that $(\Psi(t), u(t)) \in \cP$ 
and 
\[
H(\Psi(t), u(t))  = \chi(\Psi(t)) 
\sands
{\rm d} (H(\, \cdot \, , u(t))  - \chi)(\Psi(t))  = 
\]
for $0 \leq t \leq \tfin$,
\end{enumerate}
then $(\Psi, u)$ is a solution 
of the fixed time synthesis problem of duration $\tfin$ determined by $X$, $C$, and the boundary data $z_0 = \pi(\Psi(0))$ and $\zfin = \pi(\Psi(\tfin))$ .

If, in addition, $H(\Psi(0)) = 0$, then $\Psi$ determines a solution of the synthesis problem.
\end{Proposition}

\proof
For each $t \in [0, \tfin]$, synthesizability of $X$ and $C$ implies that there is
a neighborhood $\cV_t$ of $\pi(\Psi(t))$ such that  $H_t := H(\ \cdot \ , \upsilon(t)) \in \cC^1(\pi^{-1}(\cV_t))$. 
\[
0 = d(\chi - H_t)(\Psi(t))
\]
implies
\[
\dot \Psi(t) = X_\chi(\Psi(t)) = X_{H_t}(\Psi(t))
\]
and
\[
H_t(\Psi(t)) = H(\Psi(t), \upsilon(t)) = \chi(\Psi(t))
\]
for $0 \leq t \leq \tfin$. 
Hence $(\Psi, \upsilon)$ is a solution of the fixed time synthesis problem. 
If the Hamiltonian is identically zero along the trajectory, 
\prfend

\begin{Definition}
\label{mod_pot}
If \begin{enumerate}
\item
$\tilC$ is a moderation incentive,
\item
the pair $X$ and $ - \tilC$ is synthesizable,
\item
$\chi$ given by \eqr{chi_def} for  $C = - \tilC$ is $\cC^1$
\item
there is a unique map $\upsilon \in \cC^0(T^* \cS)$ such that 
\begin{enumerate}
\item
$\mbox{\em graph}(\upsilon) \subseteq \cP$, 
\item
$H(\ \cdot \ , \upsilon(\ \cdot \ )) = \chi$, and
\item
for every $\ps \in T^* \cS$, $\ps$ is a critical point of $\chi - H(\ \cdot \ , \upsilon(\ps))$,
\end{enumerate}
\end{enumerate}
then we will say that $\chi$ is a {\em moderation potential} for $\tilC$ and $X$. 
\end{Definition}

It follows immediately from Proposition \ref{build_sols} and Definition \ref{mod_pot} that 
if $\chi$ is a moderation potential for synthesizable $X$ and $\tilC$,
$\hatC \in \cC^1(\cS)$,  and $\Psi: [0, \tfin] \to T^* \cS$ is a solution of Hamilton's equations for the Hamiltonian 
$\chi  - \hatC \circ \pi$, then $(\Psi, \upsilon \circ \Psi)$ is a solution 
of the synthesis problem determined by $X$, $C = \hatC \circ \bbP_1 - \tilC$, and the boundary data
$z_0 = \pi(\Psi(0))$ and $\zfin = \pi(\Psi(\tfin))$ .
If, in addition, $H(\Psi(0)) = \hatC(\pi(\Psi(0)))$, then $\Psi$ determines a solution of the synthesis problem.

 We now show that moderation potentials exist for the family of moderation incentives constructed in 
Theorem \ref{emi}. For some subfamilies, the moderation potentials have particularly simple expressions.

\begin{Theorem}
The moderation incentives \eqr{ellip_mod_incentive} have moderation potentials 
\[
\chi_{\alpha, p}(\ps; \mu) := a_0(\ps) +  \mu(z) \, \widehat \chi_{\alpha, p}(\ell_\mu(\ps)),
\]
where $a_0(\ps) := \ps \cdot f(z)$ is the contribution of the drift field, 
\beq{mpalpha}
\widehat \chi_{\alpha, p}(r) := \left . r  \, s \lp 1 +  \smallfrac 1 {\alpha \, p} \lp r \, s^{1 - \alpha \, p} \rp^{\frac 1  {\alpha - 1}} \rp \right |_{s = \rho_{\alpha, p}^{-1}(r)}
\eeq
if \, $0 < \alpha < 1$, and 
\beq{QCC_cp}
\widehat \chi_{1, p}(r) :=  \lcb \begin{array}{ll}
\smallfrac 1 p + \smallfrac 1 q \,  r^ q  \qquad & r < 1 \medskip \\
r & r \geq 1\end{array} \right . , 
\qquad \mbox{where} \quad \smallfrac 1 p + \smallfrac 1 q = 1.
\eeq
\end{Theorem}
\prf
Setting $r = \ell_\mu(\ps)$ for notational simplicity,
Theorem \ref{emi} implies that
\beq{tmp_exp}
\frac {\chi_{\alpha, p}(\ps) - a_0(\ps)}{\mu(z)}
= \left . s \, r + \frac 1 {\alpha \, p}   \lp1 - s^p\rp^\alpha \right |_{s = \rho_{\alpha, p}^{-1}(r)}. 
\eeq
For $0 < \alpha < 1$, we can simplify this expression as follows: 
\[
r = \rho_{\alpha, p}(s) = s^{p -1 } \lp 1 - s^p \rp^{\alpha - 1}
\]
implies that 
\[
 \lp 1 - s^p \rp^\alpha  = \lp r \,  s^{p -1 } \rp^{\frac \alpha {\alpha - 1}}.
 \]
 Substituting this into \eqr{tmp_exp} and regrouping terms yields \eqr{mpalpha}.

In the case $\alpha = 1$, 
\beq{tmpit}
s \mapsto r \, s + \smallfrac 1 p \lp 1 - s^p \rp
\eeq
is the restriction of a polynomial to $\lsb 0, 1 \rsb$. The maximum of the polynomial \eqr{tmpit}
occurs at $r^{\frac 1 {p - 1}}$; hence if $r < 1$, the maximum is
\[
r^{1 + \smallfrac 1 {p - 1}} + \smallfrac 1 p \lp 1 -  r^ {\frac p {p - 1}} \rp
= \smallfrac 1 p + \lp 1 - \smallfrac 1 p \rp r^ {\frac p {p - 1}}.
\]
If $1 \leq r$, the maximum occurs at one of the endpoints;  since $\frac 1 p < 1$,
in this case, \eqr{tmpit} achieves its maximum of $r$ at $1$.
\prfend

Direct utilization of conservation of the Hamiltonian $H$ 
can simplify the analysis of the synthesis problem in many situations, However, when numerically approximating solutions of 
Hamiltonian systems, discretion must be used when combining conservation laws with discretization to avoid artificial
accelerations and related errors. We now focus on explicit use of the conservation law not to recommend it as a general purpose strategy,
but to emphasize the role of the specific value of the Hamiltonian in determining solutions satisfying given boundary conditions.
The following results play a pivotal role in our analysis of the projectile problem in \S \ref{vto}.  

We can express the optimal scalings for the moderation incentives $\tilC_{\alpha, p}$ as functions $\hatsig_{\alpha, p}$ of 
\beq{ell_eqn}
\phi(z; h) := \frac {\hatC(z) + h}{\mu(z)},
\eeq
where $h$ denotes the difference of the Hamiltonian and the drift potential at $\ps$.

\begin{Proposition}
\label{cons_quan_sols}
The optimal scaling for the moderation incentive $\tilC_{\alpha, p}$ and associated Hamiltonian
(\ref{hamap}) satisfies 
\[
\sigma_{\alpha, p}(\ps; \mu) = \hatsig_{\alpha, p}(\phi(z ; H(\ps) - a_0(\ps))),
\]
where $\hatsig_{\alpha, p}: \R \to (0, 1)$ and $\tau_{\alpha, p}: (0, 1) \to \R$ are given by 
\beq{rescale_cm_spec}
\hatsig_{\alpha, p}(\phi ) 
= \lp 1 - \tau_{\alpha, p}^{-1}(\phi))\rp^{\frac 1 p}
\qquad \mbox{for} \qquad
\tau_{\alpha, p}(w) := w^{\alpha - 1} \lp 1 + \lp \smallfrac 1 {\alpha \, p} - 1\rp w \rp
\eeq
if $0 < \alpha < 1 \leq p$ and
\beq{rescale_cm_onep}
\hatsig_{1, p}(\phi) 
= \lcb \begin{array}{ll}
{\displaystyle \lp \frac {p \, \phi - 1} {p - 1}\rp^{\frac 1 p}}& \qquad \frac 1 p \leq \phi < 1 \smallskip \\
1 & \qquad \phi \geq 1
\end{array} \right . 
\eeq
if $p > 1$.
\end{Proposition}
\prf
If $0 < \alpha < 1$ and $\lambda(\ps) \neq 0$, (\ref{sigma_alpha}) and (\ref{ell_eqn}) imply that
\beq{rescale_cm}
\widehat \chi_{\alpha, p} (\rho_{\alpha, p} (\sigma_{\alpha, p}(\ps; \mu))) = \phi(z, h).
 \eeq
The composition $\widehat \chi_{\alpha, p} \circ \rho_{\alpha, p}$ can be simplified as follows. Substituting 
\[
\rho_{\alpha, p}(s) s^{1 - \alpha \, p} = s^{(1 - \alpha) p} \lp 1 - s^p \rp^{\alpha - 1} = \lp s^{-p} - 1 \rp^{\alpha - 1}
\]
into (\ref{mpalpha}) and regrouping terms yields
\beqa
\widehat \chi_{\alpha, p}(\rho_{\alpha, p}(s)) &=& 
 \rho_{\alpha, p}(s) \, s \lp 1 +  \smallfrac 1 {\alpha \, p} \lp s^{-p} - 1\rp \rp \\
 &=& \lp 1 - s^p \rp^{\alpha - 1} s^p  \lp 1 + \smallfrac 1 {\alpha \, p} \lp s^{-p} - 1 \rp \rp \\
 &=&  \lp 1 - s^p \rp^{\alpha - 1} \lp \smallfrac 1 {\alpha \, p}  + \lp 1 -  \smallfrac 1 {\alpha \, p} \rp s^p \rp \\
 &=& \tau_{\alpha, p}\lp 1 - s^p \rp.
 \eeqa
 $\tau_{\alpha, p}$ is strictly decreasing for $0 < \alpha \leq 1$, and hence is invertible.
Solving (\ref{rescale_cm}) for $\sigma_{\alpha, p}(\ps; \mu)$ yields (\ref{rescale_cm_spec}).
 
If $\alpha = 1 < p$, then (\ref{QCC_cp}) implies 
\[
\widehat \chi_{1, p}^{-1}(\phi) =  \lcb \begin{array}{ll}
\lp q \lp \phi - \smallfrac 1 p \rp \rp^ {\frac 1 q}  \qquad & \phi < 1 \medskip \\
\phi & \phi \geq 1\end{array} \right . , 
\qquad \mbox{where} \quad \smallfrac 1 p + \smallfrac 1 q = 1.
\]
Substituting $\ell_\mu(\ps) = \widehat \chi_{1, p}^{-1}(\phi)$ in (\ref{QCC_ocp}) and simplifying yields
(\ref{rescale_cm_onep}).
\prfend

In the case $\alpha = \frac 1p$, 
 the moderation potential and optimal scaling are particularly simple. We will make use of
the following expressions in \S \ref{vto} when analyzing generalizations of the projectile example from \S \ref{intro}.

\begin{Cor}
If $\alpha = \frac 1 p < 1$, then
\beq{apone}
\shp{\chi}(\ps; \mu) - a_0(\ps) = \norm{(\ell(\ps), \mu(z))}_q = \lp \ell(\ps)^q + \mu(z)^q \rp^{\frac 1 q},
\eeq
where $\frac 1 p + \frac 1 q = 1$, and
\beq{optapone}
 \shp{\widehat\sigma}(\phi)  = \lp 1 - \phi^{-q} \rp^{\frac 1 p}.
\eeq
\end{Cor}
\prf
Setting $q = \frac p {p - 1}$ and substituting $\alpha = \frac 1 p$ and \eqr{qp_case} into \eqr{mpalpha} yields
\[
\shp{\widehat\chi}(r) = \left . r  \, s \lp 1 +  r^ {-q} \rp \right |_{s = \lp 1 + r^ {-q} \rp^{- \frac 1 p}} = r \lp 1 + r^ {-q} \rp^{\frac 1 q} =  \lp  r^ q + 1\rp^{\frac 1 q}.
\]
Hence 
\[
\shp{\chi}(\ps) - a_0(\ps) = \mu(z)\,  \shp{\widehat \chi}(\ell_\mu(\ps)) = \mu(z) \lp \ell_\mu(\ps)^q + 1 \rp^{\frac 1 q} 
= \lp \ell(\ps)^q + \mu(z)^q  \rp^{\frac 1 q}.
\]

(\ref{optapone}) follows immediately from (\ref{rescale_cm_spec}).
\prfend

\section{Vertical take-off interception with controlled velocity}
\label{vto}

To illustrate some features of moderation potentials, we return to the two dimensional system with controlled velocities
briefly considered in \S \ref{intro}, generalizing the cost functions and admissible control regions. The control problem
is designed to result in an integrable Hamiltonian system---we can express the height of the projectile and the elapsed 
time as definite integrals of functions of the horizontal position. Further specialization yields situations in which these 
definite integrals have closed form expressions in terms of elliptic integrals or logarithms, facilitating comparison of solutions
with different moderation incentives, admissible control regions, and targets. In particular, we shall see that the solutions 
change in a highly nontrivial way as the level set of Hamiltonian containing the solution changes; this illustrates the crucial
difference between the optimal control problem of unspecified duration, for which solutions must lie in the zero level set, from the fixed 
time problem, for which the appropriate level set is determined in part by the time constraint.

We first briefly recap the projectile problem from \S \ref{intro} and describe the generalizations considered here.
The task is to hit a target $(\xfin, \yfin)$, starting from an unspecified position $(x_0, 0)$ on the horizontal axis to the right of the target,
with  vertical initial velocity; the state space is $\cS = I \times \R$ for a closed interval $I \subset \R^+$ of the form $[\xfin, \infty)$
or $[\xfin, x_{\mbox{max}}]$. The velocity is the control., i.e. $(\dot x, \dot y) = u$. 

The admissible control region $\cA_{(x, y)}$ associated to $(x, y) \in \cS$ is  the closed ball of radius $\rad(x)$
centered at the origin for a given function $\rad \in \cC^0(I, \R^+)$. The unmoderated position-dependent cost term $\hatC \in \cC^0(I, \R^+)$
is a function of the horizontal component of the position. The moderation incentives have the form
$\mu(x) \tilC_{\alpha, p}$ for $\mu \in \cC^0(I, \R^+)$ satisfying $\mu(x) < \hatC(x)$ for all $x \in I$, $\alpha$ and $p$ as in \S \ref{modpot}.

We identify $T^* \cS$ with a subset of $I \times \R^2$, and denote $\ps = ((x, y), \psi)$. We abuse notation, in the interest of reminding the
reader of the invariance of key constructs, and denote quantities depending on the state variables as depending on $x$, rather
than the pair $(x, y)$, and drop the base point from $\ps$. 
The quadratic form determining the admissible control regions takes the form 
\[
Q_x(u) = \rad(x)^{-2} \norm{u}^2,
\]
and hence $L_x \psi = \rad(x) \psi$ and $\lambda(\psi) = \rad(x) \psi$ is simply a rescaling of $\psi$. 
In particular, the vertical take-off condition is equivalent to the requirement that $\ps(0)$ be vertical.

Our hypotheses were chosen so as to yield an integrable system: a pair of scalar conservation laws enable us to express 
$\psi$ as a function of $x$ and thus reduce the synthesis problem to a first order ODE solvable by quadrature. In Proposition
\ref{cons_quan_sols}  
we showed how conservation of the Hamiltonian can be used to express the optimal scaling as a function of the state
variables and the value of the Hamiltonian. 
The Hamiltonian 
\beq{hamap}
H(x, \psi) = \chi_{\alpha, p}(x, \psi) - \hatC(x)
\eeq
for the projectile system is independent of $y$; hence it follows from Noether's Theorem that
the second component of $\psi$ is a constant of the motion for  the canonical Hamiltonian system determined by $H$. 
(See, e.g. \cite{Bloch, FoM}.)
We now show that the additional conserved quantity of this system can be used to determine
 the direction of the optimal control in terms of $x$ and the  value $h$ of the Hamiltonian.
 The resulting evolution equation can be solved explicitly only in very special situations, but implicit solutions expressing
$y$ and $t$ as definite integrals depending on $x$ can be formulated as follows. 

\begin{Proposition}
\label{implicit_sols}
Let $x_0 \in I$ and $h \in \R$ satisfy
$\phi(x_0; h) \in \mbox{\rm range} [\widehat \chi_{\alpha, p}]$  for $0 < \alpha <  1 \leq p$ or $\alpha = 1 < p$.
Define $\nqm(\, \cdot \, ;  h), \wqm(\, \cdot \, ;  h): [\xfin, x_0]  \to [0, \infty)$ by
\[
\nqm(x;  h) :=  \frac {\mu(x) \widehat \chi_{\alpha, p}^{-1}(\phi(x; h))}{r(x)}
\sands
\wqm(x, x_0; h) := \lp \frac {\nqm(x_0; h)}{ \nqm(x, h)}\rp^2.
\]
If $\nqm(\, \cdot \, ;  h)$ has a strict minimum at $x_0$, then 
\beq{special_eqs_y}
y_{\alpha, p}(x; x_0, h) := \int_x^{x_0} \frac {d \xi} {\sqrt{\wqm(x_0, \xi, h) - 1}}
\eeq
and
\beq{special_eqs_t}
t_{\alpha, p}(x; x_0, h) := \int_x^{x_0} \frac {d \xi} {r(\xi)  \widehat \sigma_{\alpha, p}(\xi; h) \sqrt{1 - \wqm(\xi, x_0, h)}}
\eeq
for $\xfin \leq x < x_0$ implicitly determine the state variables of a solution in $H^{-1}(h)$ of the Hamiltonian system determined by (\ref{hamap}).
\end{Proposition}
\prf
Invertibility of $\widehat \chi_{\alpha, p}$ follows from the identity
\[
\widehat \chi_{\alpha, p}(\rho_{\alpha, p}(s)) = \tau_{\alpha, p}\lp1 - s^p \rp
\]
and the invertibility of $\rho_{\alpha, p}$ and $\tau_{\alpha, p}$. Thus 
$\ell_\mu(\ps) = \widehat \chi_{\alpha, p}^{-1}(\phi(z; h))$ and hence
\beq{norm_psi}
\norm{\psi} =  \frac {\ell(\psi)}{r(x)} =  \frac {\mu(x) \ell_\mu(\psi)}{r(x)} =  \nqm(x;  h)
\eeq
along a solution $((x, y), \psi)$ of the canonical Hamiltonian system determined by $H$ lying in $H^{-1}(h)$.

The initial condition $\dot x(0) = 0$ implies that $\psi(0) = (0, \psi_2)$, and since $\psi_2 \neq 0$ is constant and $\nqm(x;  h)$ has a strict
minimum at $x_0$, it follows that $\psi$ is always nonzero and that $\psi_1$ equals zero only when $t = 0$, then
\beq{psi_calc}
\frac {\psi}{\norm{\psi}} = \frac {1}{\norm{\psi}}\lp - \sqrt{\norm{\psi}^2 - \psi_2^2}, \psi_2 \rp 
= \lp - \sqrt{1 - \frac {\norm{\psi(0)}^2}{\norm{\psi}^2}}, \frac {\norm{\psi(0)}}{\norm{\psi}} \rp.
\eeq
(The signs are determined by  the conditions $\xfin < x_0$ and $\yfin > 0$, which imply that $\dot x$ must be negative and $\dot y$ must be positive.)
 
It follows that there are functions $y_{\alpha, p}(x; x_0, h)$ and $t_{\alpha, p}(x; x_0, h) $ such that 
\[
y_{\alpha, p}'(x; x_0, h)   = 
\frac {X_{\alpha, p}(x; h)_2}{X_{\alpha, p}(x; h)_1}= - \sqrt{\frac {\wqm(x, x_0, h)}{1 - \wqm(x, x_0, h)}} = - \frac 1 {\sqrt{\wqm(x_0, x, h) - 1}},
\]
and hence (\ref{special_eqs_y}) holds. Analogously,
\[
t_{\alpha, p}'(x; x_0, h)  =  \frac 1 {X_{\alpha, p}(x; h)_1}
\]
implies (\ref{special_eqs_t}).
\prfend

\begin{figure}[t]
\begin{center}
\includegraphics[height=1.5in]{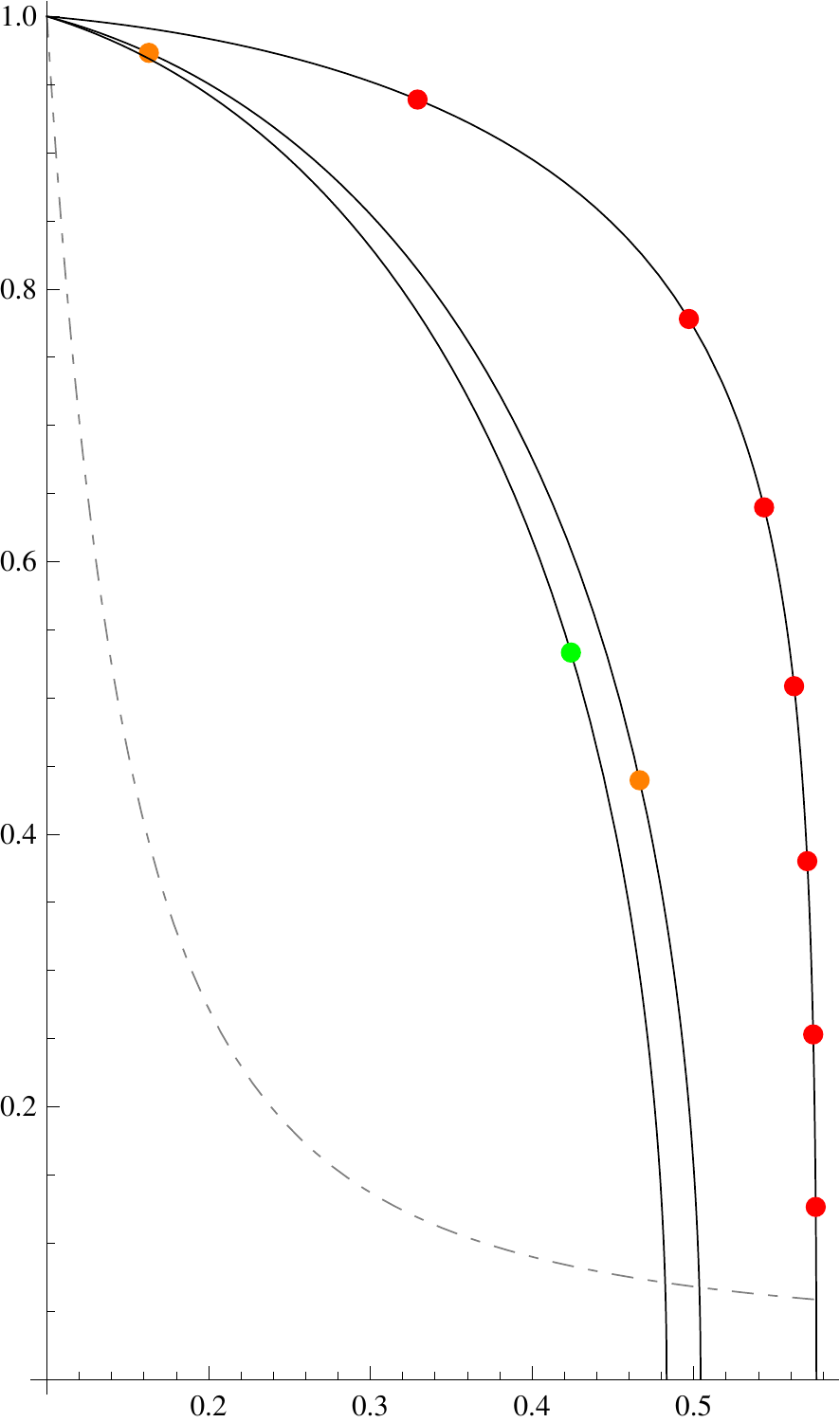} \qquad \qquad
\includegraphics[height=1.5in]{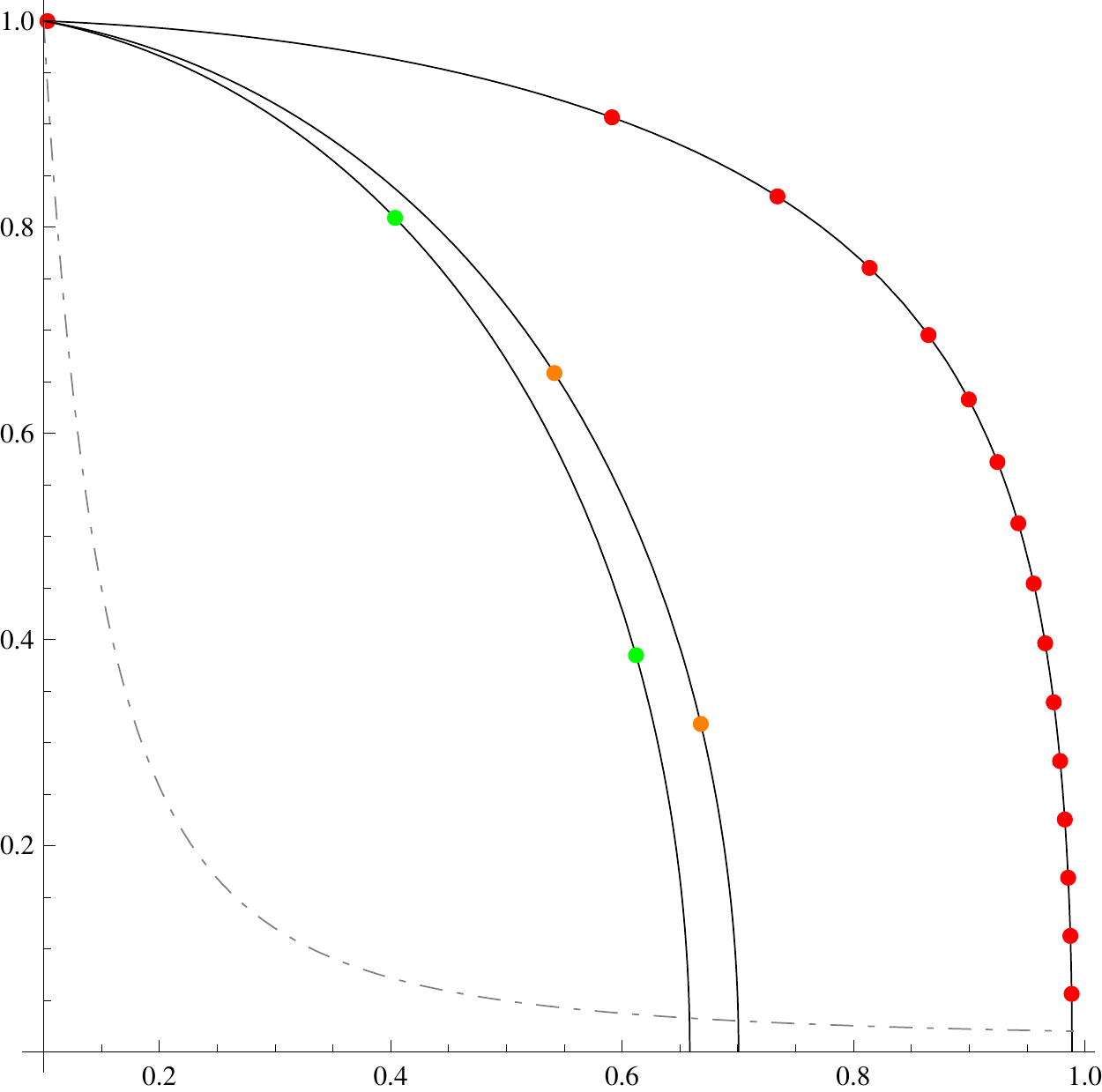} \qquad \qquad
\includegraphics[height=1.5in]{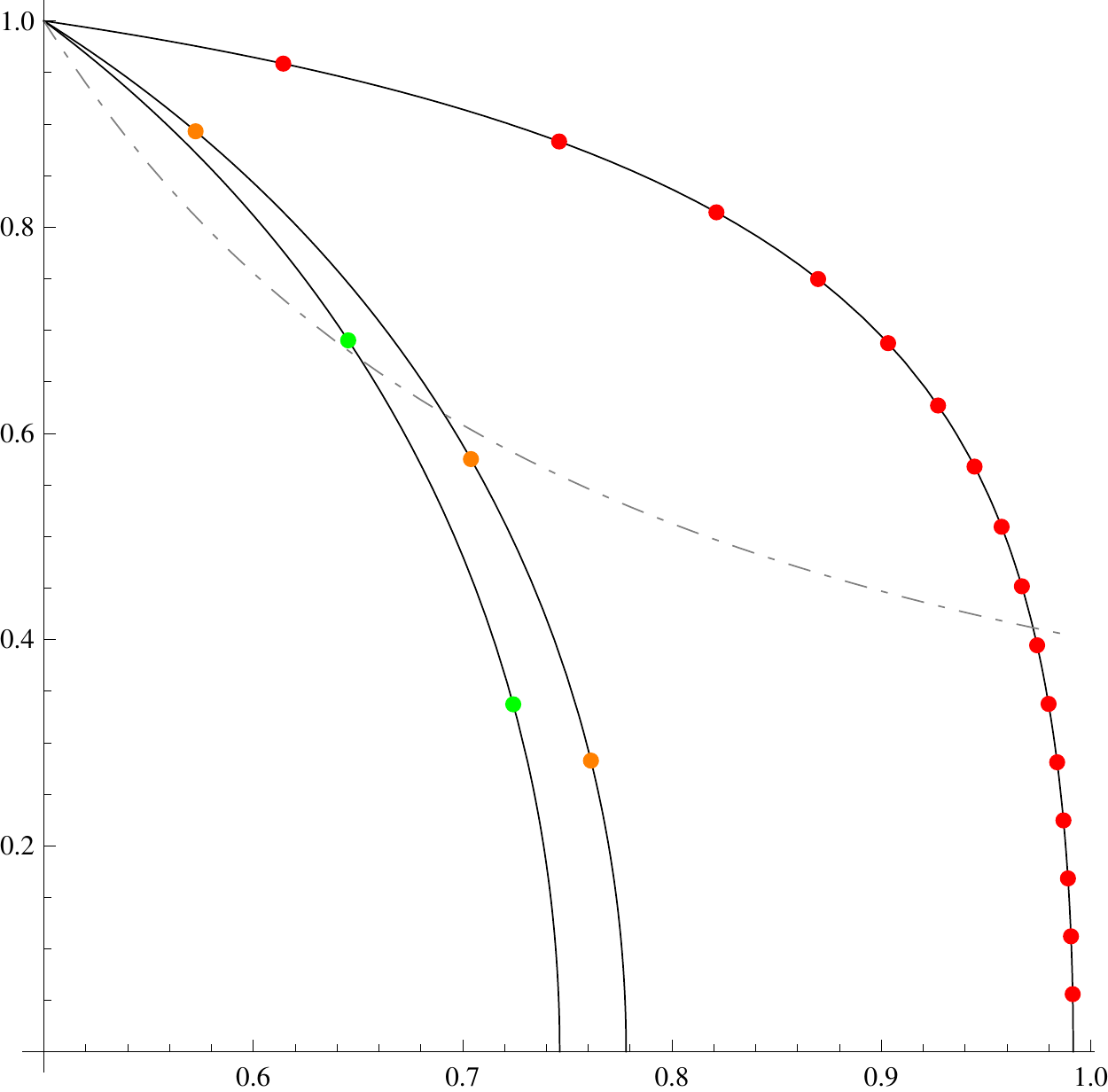} 
\end{center}
\caption{\label{log_traj_plots}Sample plots for $r (x) = \frac 1 x$. Colored dots indicate positions at $t_j := \frac j 4$. 
Left: $c = \frac 2 3$; green: $\frac \mu r = .05$, gold: $\frac \mu r = \frac 1 {\sqrt{3}}$, red: $\frac \mu r = \frac 2 {\sqrt{3}} - .05$; $\xfin = \frac 1 {10}$.
Center and right: $c = 2$; green: $\frac \mu r = .05$; gold: $\frac \mu r = 1$; red: $\frac \mu r = 1.95$;
$\xfin = \frac 1 {10}$ (center) or $\half$ (right). 
The gray dashed line indicates the normalized unmoderated cost $\frac {\hatC(x)}{\hat C(\xfin)}$. 
}
\end{figure}

\rmk
If we introduce the angle $\theta_{\alpha, p}(x; x_0, h) := \sin^{-1} \sqrt{w_{\alpha, p}(x; x_0, h)}$, $ \theta_{\alpha, p}(x; x_0, h) \in \lp \frac \pi 2, \pi \rp$,
then (\ref{psi_calc}) implies that $\psi$ has polar coordinates $(n_{\alpha, p}(x; x_0, h), \theta_{\alpha, p}(x; x_0, h))$ at $(x, y)$ and 
\[
y_{\alpha, p}'(x; x_0, h) = \tan \theta_{\alpha, p}(x; x_0, h).
\]
However, we have found it more convenient in specific calculations to work with $w_{\alpha, p}$.
\rmkend

Proposition \ref{implicit_sols} provides implicit equations for solutions of the Hamiltonian system with Hamiltonian (\ref{hamap}).
Such solutions only qualify as solutions of the synthesis problem if additional conditions on the parameters $x_0$ and
$h$ are satisfied. 
\begin{description}
 \item[\it Synthesis problem.]
 The projectile must strike the target and lie in the zero level set of the Hamiltonian. The initial position $x_0$ determines
a solution of the synthesis problem $\quad \Longleftrightarrow \quad y_{\alpha, p}(\xfin; x_0, 0) = \yfin$.
 \item[\it Fixed time synthesis problem.]
 The projectile must strike the target at the specified time $\tfin$. The initial position $x_0$ and Hamiltonian value
$h$ determine
a solution of the time $\tfin$ synthesis problem $\quad \Longleftrightarrow \quad y_{\alpha, p}(\xfin; x_0, h) = \yfin$ \ \ 
and \ \ $t_{\alpha, p}(\xfin; x_0, h) = \tfin$.
\end{description}

\rmk
Given $x_0$ and $h$ such that $y_{\alpha, p}( \, \cdot \, ; x_0, h)$  and $t_{\alpha, p}( \, \cdot \, ; x_0, h)$ are well-defined 
on $[\xfin, x_0]$ and $y_{\alpha, p}(\xfin; x_0, h) = \yfin$, one can, of course, {\it a posteriori} specify $t_{\alpha, p}(\xfin ; x_0, h)$ as
the desired duration, thereby obtaining a solution of the corresponding fixed time optimal control problem. 
However, if there is a family of pairs $(x_0, h)$ determining solutions of different durations that all strike the target
and only one of these solutions will be implemented, some criterion for selecting that solution must be established. 
\rmkend

\begin{figure}[t]
\begin{center}
\includegraphics[height=1.5in]{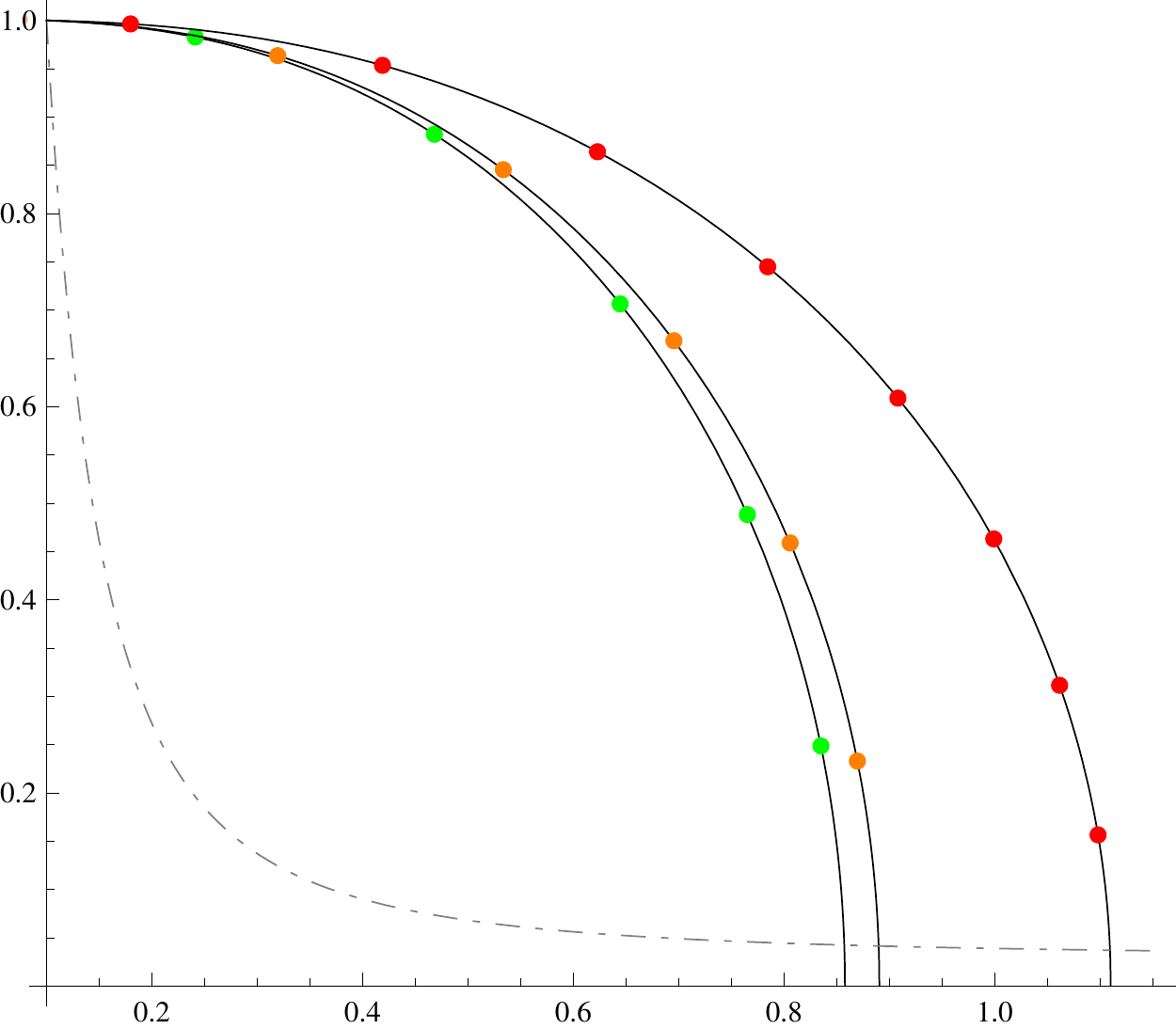} \qquad \quad 
\includegraphics[height=1.5in]{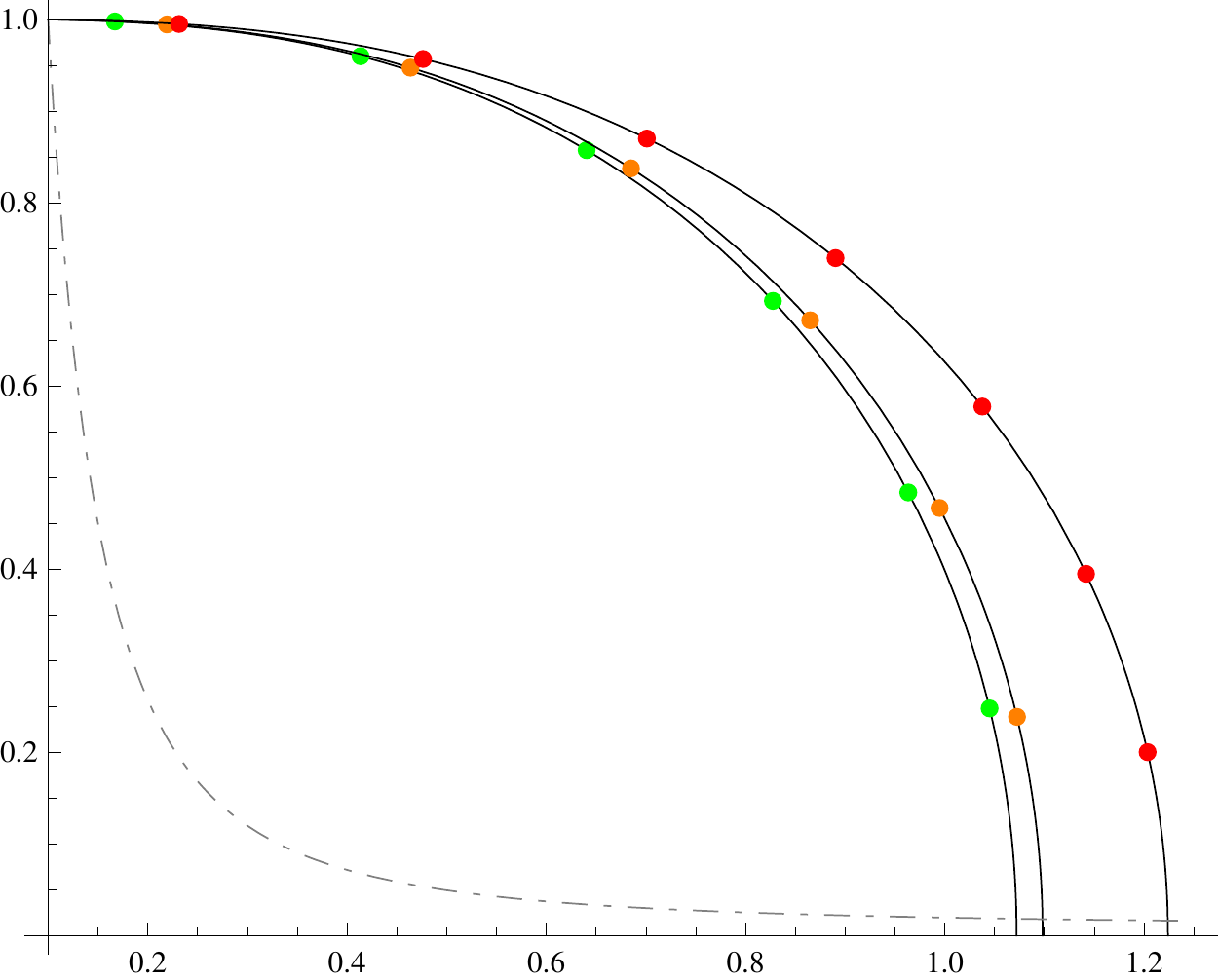} \qquad \quad
\includegraphics[height=1.5in]{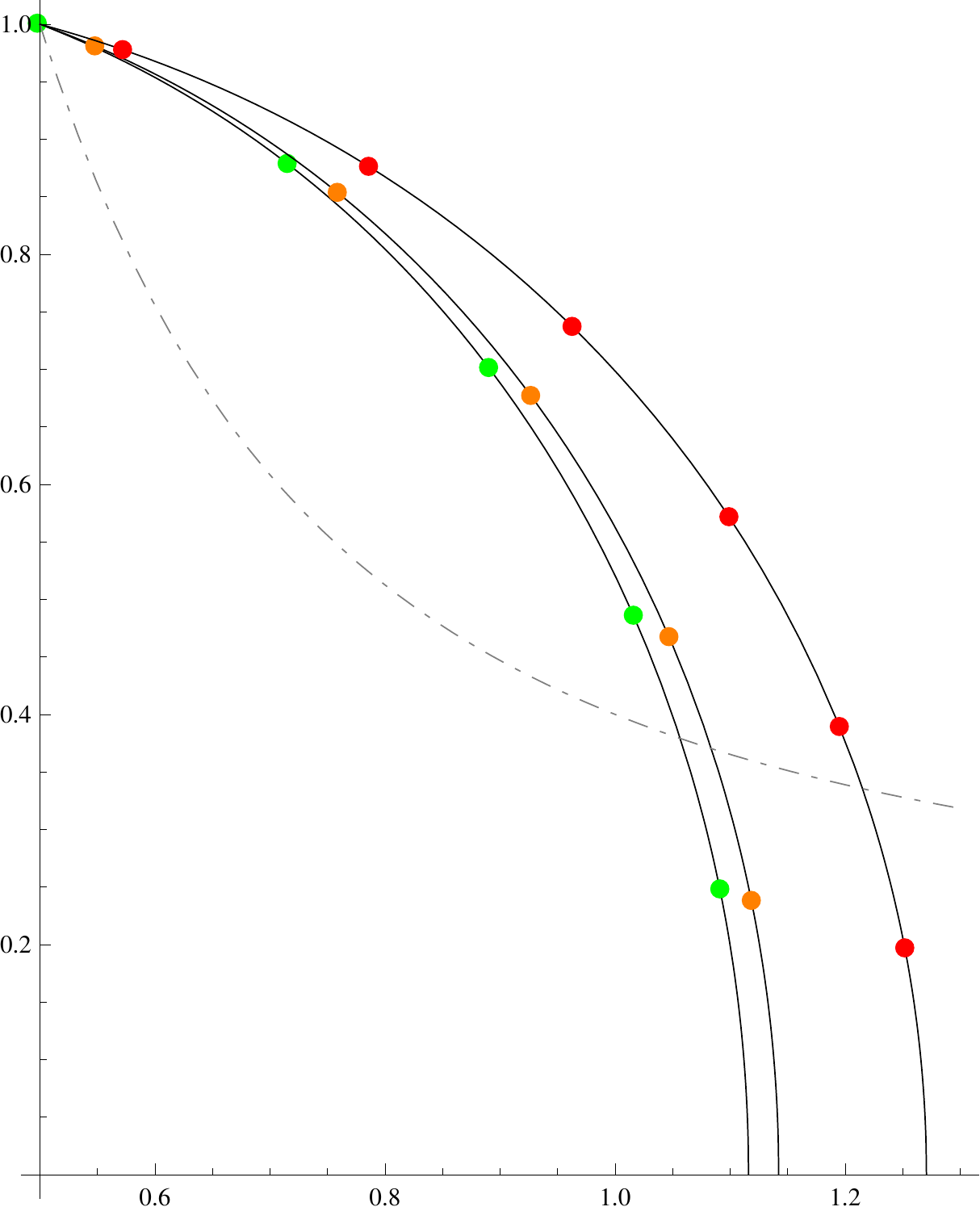} 
\end{center}
\caption{\label{ell_traj_plots}Sample trajectories for $r \equiv 1$. Colored dots indicate positions at $t_j := \frac j 4$. 
Left: $c = \frac 2 3$; green: $\frac \mu r = .05$, gold: $\frac \mu r = \frac 1 {\sqrt{3}}$, red: $\frac \mu r = \frac 2 {\sqrt{3}} - .05$; $\xfin = \frac 1 {10}$.
Center and right: $c = 2$; green: $\frac \mu r = .05$; gold: $\frac \mu r = 1$; red: $\frac \mu r = 1.95$;
$\xfin = \frac 1 {10}$ (center) or $\half$ (right). The gray dashed line indicates the normalized unmoderated cost $\frac {\hatC(x)}{\hat C(\xfin)}$. 
}
\end{figure}

\subsection{$\alpha = \half$, $p = 2$, and constant $\frac \mu r$}

The expressions for $y_{\alpha, p}$ and $t_{\alpha, p}$ as functions of $x$ take a particularly simple form if $\alpha = \half$, $p = 2$,
 and $\mu$ is a positive rescaling of $\rad$. 

\begin{Cor}
\label{mr2}
If  $\phi(x; h) < 1$ for $\xfin \leq x < x_0$, define  
 \beq{special_eqs_omega}
v(x_0, h) := \sqrt{1 - \phi(x_0; h)^{-2}}
 \eeq
 and
  \beq{special_eqs_tily_2mr}
\tilde y(x; x_0, h) := \int_x^{x_0}\frac {d\xi}{\sqrt{\lp \frac {\phi(\xi; h)}{\phi(x_0; h)} \rp^2 - 1}} 
\eeq
for $\xfin \leq x < x_0$. Then
\beq{special_eqs_yt_2mr}
\sht{y}(x; x_0, h) =  v(x_0, h) \, \tilde y(x ; x_0, h)
\eeq
and
\beq{special_eqs_t_2mr}
\sht{t}(x; x_0, h) = \int_x^{x_0} \frac {d \xi} {\rad(u) \sqrt{1 - \lp \frac {\phi(x_0; h)}{\phi(\xi; h)} \rp^2}}.
\eeq
 \end{Cor}

\prf
${\displaystyle \widehat \chi_{\half, 2}^{-1}(\phi) = \sqrt{1 - \phi^2}}$ and $\frac \mu r = \mbox{constant}$ imply that
\[
\sht{w}(x_0, x; h)  = \lp \frac {\widehat \chi_{\half, 2}^{-1}(\phi(x; h))}{\widehat \chi_{\half, 2}^{-1}(\phi(x_0; h))} \rp^2 
= \frac {1 - \phi(x; h)^2}{1 - \phi(x_0; h)^2} 
= \frac {\displaystyle \lp \frac {\phi(x; h)} {\phi(x_0; h)} \rp^2 - 1}{1 - \phi(x_0; h)^{-2}} + 1.
\] 
Analogously, 
\[
1- \sht{w}(x,x_0;  h) 
= \frac {\displaystyle 1 - \lp \frac {\phi(x_0; h)}{\phi(x; h)} \rp^2 }{\displaystyle 1 - \phi(x; h)^{-2}} 
= \frac {\displaystyle 1 - \lp \frac {\phi(x_0; h)}{\phi(x; h)} \rp^2 }{\sht{\widehat \sigma}(x; h)^2}. 
\]
Substituting these expressions and  $p = q = 2$
into (\ref{special_eqs_y}) and (\ref{special_eqs_t}) yields  (\ref{special_eqs_yt_2mr}).

The vertical component of the velocity at position $(x, y(x))$ is given by
\[
\frac {\sht{y}'(x)}{\sht{t'}(x)} = v(x_0, h) \rad(x) \sqrt{\frac  {1 - \lp \frac {\phi(x_0; h)}{\phi(x; h)}\rp^2} {\lp \frac {\phi(x; h)}{\phi(x_0; h)}\rp^2 - 1}}
= v(x_0, h) \rad(x)\frac {\phi(x_0; h)}{\phi(x; h)}.
\]
\prfend

Corollary \ref{mr2} reveals several distinctive features of this special situation.
\begin{itemize}
\item
The relationship between $t$ and $x$ does not depend on the value of the constant ratio $\frac \mu \rad$.
\item
The function $\tilde y$ determines the optimal solution of the unmoderated problem, corresponding to $\mu \equiv 0$,
with initial velocity on the boundary of the admissible control region.
\item
The relationship between $y$ and $x$ depends on the ratio $\frac \mu \rad$ only via the scaling factor $v(x_0, h)$.
\item
$x_0$ determines a trajectory with energy $h$ passing through the point $(\xfin, \yfin)$  iff $(1, \yfin)$ lies on the positive quadrant of the ellipsoid 
with principal axes $\phi(x_0; h)$ and $\tilde y(\xfin; x_0, h)$. 
\end{itemize}

If $\phi(x_0; 0)$ and $\tilde y(\xfin; x_0, 0)$ determine a family of non-intersecting ellipsoids parametrized by $x_0$, then the synthesis 
problem with target $(\xfin, \yfin)$ has a unique solution for each admissible value of the ratio $\frac \mu r$, with initial position $(x_0, 0)$
and initial velocity $(0, v(x_0, 0) \rho(x_0))$, for the unique value of $x_0$ such that the ellipse with principal axes
$\phi(x_0; 0)$ and $\tilde y(\xfin; x_0, 0)$ passes through $(1, \yfin)$.

We now further specialize, considering the position-dependent cost term $\hatC(x) = \frac c{2 \, x^2} + 1$ and 
admissible control region radius functions $\rad(x) = \frac 1 x$ or $\rad(x) = 1$. In these cases
we can explicitly express $\tilde y$ and $t$ as functions of $x$ in terms of logarithms (for $\rho(x) = \frac 1 x$) or elliptic integrals (for $\rho(x) \equiv 1$).
We present the solutions only for $h = 0$, corresponding to solutions of the synthesis problem; 
the expressions for nonzero $h$ are similar, but involve somewhat messier coefficients. 

If we set ${\eta(x, x_0) := \sqrt{\lp \frac c {2 \, x_0} \rp^2 - x^2}}$, then 
\[
\tilde y(x; x_0, 0) =  \lp \frac c {2 \, x_0} + x_0 \rp \lp \ln \lp \eta(x, x_0) + \sqrt{x_0^2 - x^2} \rp -  \ln \eta(x_0, x_0) \rp
\]
and
\[
2 \, \sht{t}(x; x_0, 0) =  \lp \frac c {2 \, x_0} + x_0 \rp \tilde y(x; x_0, 0) - \eta(x, x_0) \sqrt{x_0^2 - x^2}
\]
if $\rad(x) = \frac 1 x$.

Trajectories for some representative values of $c$, $\frac \mu r$, $\xfin$, and $x_0$, with $\rad(x) = \frac 1 x$ and $\yfin = 1$, are shown in Figure \ref{log_traj_plots}. 
Note that the more moderate the strategy, the further $x_0$ is from the target and the slower the initial ascent. Trajectories with moderation
factor near the maximum allowable value show a slow, nearly vertical early phase, executed in relatively `safe' territory (i.e. relatively small values
of $\hatC(x)$) followed by a rapid, nearly horizontal late phase; those with low moderation factor launch closer to the target and rapidly pursue
a more rounded path. The moderation factor does not correspond to increased or reduced sensitivity to risk, but influences the approach to reducing 
risk---the more moderate solution takes more time and travels a longer path overall, but in doing so, is able to devote most of its (constrained) speed 
to nearly horizontal motion when moving through the high-risk zone near the target.The differences as the moderation factor are smaller
if the risk is lower, either due to a smaller value of the risk factor $c$ or to relatively large $\xfin$, and hence relatively small variation in risk
from $x_0$ to $\xfin$. Finally, note that for this specific system, changes in the moderation factor $\frac \mu r$  result in relatively small changes in the optimal 
trajectory until $\frac \mu r$ is close to the maximum value. 

 \bigskip

If we let ${\mathcal E}_E$ and ${\mathcal E}_F$ denote the incomplete elliptic integrals of the first and second kind, and define
\[
\tilde \gamma_\pm(u; k) := \cE_F(\sin^{-1} u; k) \pm  \cE_E(\sin^{-1} u;  k),
\qquad
k(x_0) := - \lp 1 + \frac {4 x_0^2} c \rp,
\]
and
\[
\gamma_\pm(x; x_0) := x_0 \lp \tilde \gamma_\pm \lp \smallfrac x {x_0}; k(x_0) \rp -  \tilde \gamma_\pm(1; k(x_0)) \rp, 
\]
then 
\[
\tilde y(x; x_0, 0) = \frac {\gamma_-(x; x_0)}{ 1 + \frac 1 {\hatC(x_0)}} 
\sands
2 \, \sht{t}(x; x_0, 0) =  \gamma_+ (x; x_0) + \smallfrac 1 {k(x_0)} \gamma_- (x; x_0)
\]
if $\rad \equiv 1$.

Trajectories for some representative values of $c$, $\frac \mu r$, $\xfin$, and $x_0$, with $\rad(x) \equiv 1$ and $\yfin = 1$, are shown in Figure \ref{ell_traj_plots}. 
As before, the more moderate the strategy, the further $x_0$ is from the target and the slower the ascent. However, since the admissible
control region is the unit ball for all values of $x$, there is much less dramatic variation in the speed along any given solution and in the
paths of the different solutions. All of the trajectories trace follow paths that are nearly, but not exactly, elliptical. 

\subsection{$\alpha = 1$, $p = 2$, and constant $\frac \mu r$}

We now consider the parameters values and functions used in \S \ref{warm_up}: $C_{\rm mi} = \tilde C_{1, 2}( \, \cdot \, , \mu) - \hatC$
for constant $\mu$ and $\hatC$ as above and $C_{\rm ke} =  C_{\rm mi} - 1 + \frac \mu 2$. Thus solutions of the synthesis problem
for $C_{\rm ke}$ correspond to solutions of an appropriate fixed time synthesis problem for $C_{\rm mi}$.
We derive the solutions of the synthesis problem for a general $\hatC$ depending only on $x$ and Hamiltonian value $h$ before
specializing to $C_{\rm mi}$ and $C_{\rm ke}$.

When $\alpha = 1$ and $p = 2$, the condition that the instantaneous cost be positive everywhere imposes the inequality
$\phi(x; 0) > \half$, and (\ref{rescale_cm_onep}) takes the form
\beq{params12}
\widehat \chi_{1, 2}^{-1}(\phi)
= \lcb \begin{array}{cl}
\sqrt{2 \,  \phi - 1}\quad & \mbox{if} \quad \half \leq \phi < 1 \smallskip \\
\phi &  \mbox{if} \quad \phi \geq 1
\end{array} \right . 
\eeq
and
\[
\widehat \sigma_{1, 2}(x; h) 
= \min \lcb \sqrt{2 \,  \phi(x; h) - 1}, 1 \rcb.
\]

For simplicity, we consider only trajectories such that either $\half \leq \phi(x; h) \leq 1$ for $\xfin \leq x \leq x_0$ or
$1 \leq \phi(x; h)$ for $\xfin \leq x \leq x_0$; determining more general solutions involves patching together solutions
of these kinds. (\ref{params12}) implies that
\beq{w12}
w_{1, 2}(x; x_0, h) = \lcb \begin{array}{ll}
{\displaystyle \frac {2 \, \phi(x_0; h) - 1}{2 \, \phi(x; h) - 1} }& \half \leq \phi(x_0; h) < \phi(x; h) \leq 1 \smallskip \\
{\displaystyle\lp \frac {\phi(x_0; h)}{\phi(x; h)} \rp^2 }& 1 \leq  \phi(x_0; h) 
\end{array} \right . 
\eeq
for constant $\frac \mu r$.

If $1 \leq \phi(x_0; h)$, comparing (\ref{special_eqs_y}) to (\ref{special_eqs_tily_2mr}) and (\ref{special_eqs_t})
to (\ref{special_eqs_yt_2mr}) shows that in this situation
 \[
 y_{1, 2}(x; x_0, h) = \tilde y(x; x_0, h) = \frac { y_{{\tiny 1 \over 2}, 2}(x; x_0, h)}{v(x_0; h)} 
 \sands
t_{1, 2}(x; x_0, h) = \sht{t}(x; x_0, h).
\]

If $\half \leq \phi(x; h) \leq 1$ for $\xfin \leq x \leq x_0$, then 
\beq{slow12}
\frac {y_{1, 2}(x; x_0, h)}{\widehat \sigma_{1, 2}(x_0; h)}
= t_{1, 2}(x; x_0, h) =   \int_x^{x_0} \frac {d \xi}{\sqrt{2 (\phi(\xi; h) - \phi(x_0; h))}}.
\eeq
If we further specialize to the case $\hatC(x) =  1 + \frac c {2 \, x^2}$, $r \equiv 1$, then
\[
t_{1, 2}(x; x_0, h) = \sqrt{\frac {x_0^2 - x^2} {\zeta(x_0)}}
\sands
\widehat \sigma_{1, 2}(x_0; h) = \zeta(x_0) + \smallfrac {2 (1 + h)} \mu - 1
\qquad \mbox{for} \qquad \zeta(x_0) :=  \frac {c}{\mu \, x_0^2},
\]
where $\mu$ denotes the constant value of the moderation factor.  It follows that
the projectile paths are segments of ellipses centered at the origin.
We can easily express $z = (x, y)$ explicitly as a function of $t$ in this case: 
\[
z_{1, 2}(t; h) = \lp \sqrt{x_0^2 - \zeta(x_0) t^2},  \widehat \sigma_{1, 2}(x_0; h) t \rp.
\]
Setting $h = 0$ yields the state information of the synthesis problem for $C_{\rm mi}$, while setting $h = \frac \mu 2 - 1$ gives the corresponding
information for $C_{\rm ke}$.

\rmk
The graphs of $y_{1, 2}(\, \cdot \, ; x_0, h)$ are very nearly elliptical if $1 \leq \phi(x_0; h)$, but do not exactly coincide with segments of ellipses. 
\rmkend

\bibliographystyle{plain}
\bibliography{control}

\begin{thebibliography}{10}

\bibitem{FoM}
R.~Abraham and J.E. Marsden.
\newblock {\em Foundations of mechanics}.
\newblock Benjamin/Cummings Publishing Company, Reading, MA, 1978.

\bibitem{BW}
J.~Baillieul and J.C. Willems.
\newblock {\em Mathematical Control Theory}.
\newblock Springer, 1999.

\bibitem{BB}
Enrico Bertolazzi, Francesco Biral, and Mauro~Da Lio.
\newblock Real-time motion planning for multibody systems: real life
  application examples.
\newblock {\em Multibody Syst. Dyn.}, 17:119--139, 2007.

\bibitem{Bloch}
Anthony Bloch.
\newblock {\em Nonholonomic Mechanics and Control}.
\newblock Springer, 2003.

\bibitem{BG}
J.~F. Bonnans and Th. Guilbaud.
\newblock Using logarithmic penalties in the shooting algorithm for optimal
  control problems.
\newblock {\em Optimal Control Appl. Methods}, 24:257--278, 2003.

\bibitem{life}
Ralph Crane~(photographer).
\newblock A copycat astronaut.
\newblock {\em L{I}{F}{E}}, August 16:77--78, 1968.

\bibitem{KS}
T.R. Kane and M.P. Scher.
\newblock A dynamical explanation of the falling cat phenomenon.
\newblock {\em Int. J. Solids Structures}, 5:663--670, 1969.

\bibitem{Leimkuhler}
Benjamin Leimkuhler and Sebastian Reich.
\newblock {\em Simulating Hamiltonian Dynamics}.
\newblock Cambridge University Press, 2005.

\bibitem{Lewis_cons}
Debra Lewis.
\newblock Optimal control with moderation incentives.
\newblock Preprint, arXiv:1001.0211v1 [math.OC], 2010.

\bibitem{LO1}
Debra Lewis and Peter Olver.
\newblock Geometric integration algorithms on homogeneous manifolds.
\newblock {\em Foundations of Computational Mathematics}, 2:363--392, 2002.

\bibitem{LS}
Debra Lewis and J.C. Simo.
\newblock Conserving algorithms for the dynamics of {H}amiltonian systems on
  {L}ie groups.
\newblock {\em J. Nonlin. Sci.}, 4:253--299, 1994.

\bibitem{Marey}
E.-J. Marey.
\newblock M\'echanique animale.
\newblock {\em La Nature}, 1119:569--570, 1894.

\bibitem{Mont92}
Richard Montgomery.
\newblock Optimal control of deformable bodies and its relation to gauge
  theory.
\newblock {\em Math. Sci. Res. Inst. Publ.}, 22:403--438, 1991.

\bibitem{Mont93}
Richard Montgomery.
\newblock Gauge theory of the falling cat.
\newblock {\em Fields Institute Communications}, 1:193--218, 1993.

\bibitem{NVdS}
H.~Nijmeijer and A.J. van~der Schaft.
\newblock {\em Nonlinear Dynamical Control Systems}.
\newblock Springer-Verlag, 1990.

\bibitem{PBGM}
L.S. Pontryagin, V.G. Boltyanskii, R.V. Gamkrelidze, and E.F. Mishchenko.
\newblock {\em The Mathematical Theory of Optimal Processes}.
\newblock Interscience Publishers, 1962.

\bibitem{Sontaga}
Eduardo~D. Sontag.
\newblock Integrability of certain distributions associated with actions on
  manifolds and applications to control problems.
\newblock In {\em Nonlinear controlability and optimal control}, pages 81--131.
  Marcel Dekker, Inc., 1990.

\bibitem{Sontagb}
Eduardo~D. Sontag.
\newblock {\em Mathematical Control Theory: Deterministic Finite Dimensional
  Systems}.
\newblock Springer, 1998.

\bibitem{synode}
{SYNODE} publication database. {\tt
  http://www.math.ntnu.no/num/synode/bibliography.php}, 2003.

\end{thebibliography}
\end{document}